\documentclass[graybox]{svmult}

\usepackage{amsfonts}
\usepackage{amssymb}
\usepackage{amsmath}
\usepackage{eucal}

\newcommand{\ellipK}{\mathcal{K}}
\newcommand{\ellipE}{\mathcal{E}}
\newcommand{\arth}{\operatorname{arth}}
\newcommand{\R}{\mathbb{R}}

\usepackage{mathptmx}       
\usepackage{helvet}         
\usepackage{courier}        
\usepackage{type1cm}        
\usepackage{makeidx}         
\usepackage{graphicx}        
\usepackage{multicol}        
\usepackage[bottom]{footmisc}

\makeindex             

\begin{document}

\title*{Topics in special functions III}
\author{Glen D. Anderson, Matti Vuorinen, and Xiaohui Zhang }
\institute{Glen D. Anderson \at Department of Mathematics, Michigan State University, East Lansing, MI 48824, USA,
            \email{anderson@math.msu.edu}
            \and Matti Vuorinen \at Department of Mathematics and Statistics, University of Turku, Assistentinkatu 7,      FIN-20014, Finland,  \email{vuorinen@utu.fi}
            \and Xiaohui Zhang \at Department of Mathematics and Statistics, University of Turku, Assistentinkatu 7,      FIN-20014, Finland,  \email{xiazha@utu.fi}}
%
%
\maketitle

\abstract*{The authors survey recent results in special functions of classical analysis and geometric function theory, in particular the circular and hyperbolic functions, the gamma function, the elliptic integrals, the Gaussian hypergeometric function, power series, and mean values.}

\section{Introduction}
\label{sec:1}

The study of quasiconformal maps led the first two authors in their joint work with
M. K. Vamanamurthy to formulate open problems or questions involving special functions \cite{avv6,avv9}.  During the past two decades, many authors have contributed to the solution of these problems.
However, most of the problems posed in \cite{avv6} are still open.

The present paper is the third in a series of surveys by the first two
authors, the previous papers \cite{avv3,avv5} being written jointly with the
late M. K. Vamanamurthy. The aim of this series of surveys is to review the
results motivated by the problems in \cite{avv6,avv9} and related
developments during the past two decades.
In the first of these we studied classical special functions, and in the next we focused on special functions occurring in the distortion theory of quasiconformal maps.  Regretfully, Vamanamurthy passed away in 2009, and the remaining authors dedicate the present work as a tribute to his memory. For an update to the bibliographies of
\cite{avv3} and \cite{avv5} the reader is referred to \cite{av1}.

In 1993 the following monotone rule was derived \cite[Lemma 2.2]{avv1}.  Though simple to state and easy to prove by means of the Cauchy Mean Value Theorem, this l'H\^opital Monotone Rule (LMR) has had wide application to special functions by many authors.  Vamanamurthy was especially skillful in the application of this rule. We here quote the rule as it was restated in \cite[Theorem 2]{avv4}.

\begin{theorem}\label{thm:LMR}
\rm{(l'H\^opital Monotone Rule).} Let  $-\infty < a < b < \infty$, and let  $f,g: [a,b] \to \mathbb{R}$ be continuous functions that are differentiable on  $(a,b)$, with $f(a)=g(a)=0$ or $f(b)=g(b)=0.$  Assume that  $g'(x)\ne 0$ for each  $x\in (a,b).$  If  $f'/g'$  is increasing (decreasing) on  $(a,b)$, then so is $f/g$.
\end{theorem}

Theorem \ref{thm:LMR} assumes that
$a$ and $b$ are finite, but the rule can be extended easily by similar methods to the case
where $a$ or $b$ is infinite. The l'H\^opital Monotone Rule has been used effectively in the study of the monotonicity of a quotient of two funcions. For instance, Pinelis' note \cite{pi2} shows the potential of the LMR. As a complement to Pinelis' note, the paper \cite{avv4} contains many applications of LMR in calculus.
Also the history of LMR is reviewed there.

In this survey we give an account of the work in the special functions of classical analysis and geometric function theory since our second survey.  In many of these results the  LMR was an essential tool.  Because of practical constraints, we have had to exclude many fine papers and have limited our bibliography to those papers most closely connected to our work.

The aim of our work on special functions has been to solve open problems in quasiconformal mapping theory.
In particular, we tried to settle Mori's conjecture  for quasiconformal mappings \cite{mo} (see also \cite[p. 68]{lv}). For the formulation of this problem, let $K>1$ be fixed and let $M(K)$ be the least constant such that
$$
|f(z_1)-f(z_2)|\le M(K)|z_1-z_2|^{1/K},\,\ \ {\rm for\ all}\ z_1,z_2\in B,
$$
for every $K$-quasiconformal mapping $f:B \to B$ of the unit disk $B$ onto itself with $f(0)=0$.
A. Mori conjectured in 1956 that $M(K) \le 16^{1-1/K}\,.$ This conjecture is still open in 2012. Some of the open problems that we found will be discussed in the last section.

\section{Generalizations of Jordan's inequality}

The LMR application list, begun in \cite{avv4}, led to the Master's thesis of M. Visuri, on which \cite{kvv} is based.
Furthermore, applications of LMR to trigonometric inequalities were given in \cite{kvv}. Numerous further applications to
trigonometric functions were found by many authors, and some of these papers are reviewed in this section and the next.

By elementary geometric methods one can prove that
$$
\frac{2}{\pi}\le \frac{\sin x}{x}<1,\ \ 0<x\le\frac{\pi}{2},
$$
a result known as \emph{Jordan's inequality}.
In a recent work, R. Kl\'en et al \cite{kvv} have obtained the inequalities
$$
\cos^2\frac{x}{2}\le\frac{\sin x}{x}\le \cos^3\frac{x}{3},\ \ x\in (-\sqrt{27/5},\sqrt{27/5}),
$$
and
$$
 {\rm cosh^{1/4}\  x} < \frac{{\rm sinh}\ x}{x}<{\rm cosh^{1/2}}x,\ \ x\in (0,1).
$$
Inspired by these results,  Y.-P. Lv, G.-D. Wang, and Y.-M. Chu \cite{lwc}  proved that, for  $a = (\log(\pi/2))/\log\sqrt{2}\approx1.30299$,
$$
\cos^{4/3}\frac{x}{2}<\frac{\sin x}{x}<\cos^a\frac{x}{2}, \ \ x\in(0,\pi /2),
$$
where $4/3$ and  $a$ are best constants, and that for $b=(\log\rm{sinh} 1)/(\log\rm{cosh} 1)\approx 0.372168,$
$$
 {\rm cosh^{1/3}}  x < \frac{{\rm sinh}\ x}{x}<{\rm cosh}^{b}x,\ \ x\in (0,1),
$$
where  $1/3$ and $b$ are best constants.

Many authors have generalized or sharpened Jordan's inequality, either by replacing the bounds by finite series or hyperbolic functions, or by obtaining analogous results for other functions such as hyperbolic or Bessel functions.  The comprehensive survey paper by F. Qi et al \cite{qng} gives a clear picture of these developments as of 2009.   For example, in 2008 D.-W. Niu et al \cite{nhcq} obtained the sharp inequality
$$
\frac{2}{\pi}+\sum_{k=1}^{n}\alpha_k(\pi^2-4x^2)^k\le\frac{\sin x}{x}\le\frac{2}{\pi} +\sum_{k=1}^{n}\beta_k(\pi^2-4x^2)^k,\, \, 0<x\leq{\pi}/{2},
$$
for each natural number $n$, with best possible constants $\alpha_k$ and $\beta_k$.  That same year S.-H. Wu and H. M. Srivastava \cite{ws2} obtained upper and lower estimates on $(0,\pi/2]$ for $(\sin x)/x$ that are finite series in powers of  $(x-\theta),$ where $\theta\in [x,\pi/2],$ while L. Zhu \cite{z3} obtained bounds as finite series in powers of $(\pi^2-4x^2)$.  L. Zhu \cite{z2}  obtained bounds for $(\sin x)/x$  as finite series in powers of $(r^2-x^2)$ for  $0 < x \le r \le \pi/2$, yielding a new infinite series
$$
\frac{\sin x}{x} = \sum_{n=0}^{\infty}a_n(r^2-x^2)^n, \ \ {\rm for}\ \ 0<|x|\le r\le \pi/2.
$$
S.-J. Yang \cite{YangSh} showed that a function $f$ admits an infinite series expansions of the above type if and only if $f$ is analytic and even.

In 2011 Z.-H. Huo et al \cite{hncq} obtained the following generalization of Jordan's inequalities:
$$
\sum_{k=1}^n\mu_k(\theta^t-x^t)^k\le \frac{\sin x}{x} - \frac{\sin \theta}{\theta}\le \sum_{k=1}^n\omega_k(\theta^t-x^t)^k
$$
for $t\ge 2$, $n\in \mathbb{N}$, and  $0<x\leq\theta<\pi$, where the coefficients $\mu_k$ and $\omega_k$ are defined recursively and are best possible.

More recently, in 2012, C.-P. Chen and L. Debnath \cite{cd} have proved that, for $0<x\le\pi/2,$
$$
f_1(x)\le \frac{\sin x}{x} \le f_2(x),
$$
where
$$
f_1(x) = \frac{2}{\pi} + \frac{2 \pi^{-\theta-1}}{\theta} (\pi^{\theta}-(2x)^{\theta})+\frac{(-\pi^2+4+4\theta)\pi^{-2\theta-1}}{4\theta^2}(\pi^{\theta}-(2x)^{\theta})^2
$$
and
$$
f_2(x)=\frac{2}{\pi} + \frac{2 \pi^{-\theta-1}}{\theta} (\pi^{\theta}-(2x)^{\theta})+\frac{((\pi-2)\theta-2)\pi^{-2\theta-1}}{\theta}(\pi^{\theta}-(2x)^{\theta})^2,
$$
for any $\theta\ge 2$, with equality when $x=\pi/2$.

In a recent work J. S\'andor \cite{s1} (see also \cite[Paper 9]{s4})  proved that $h(x)\equiv[\log(x/\sin x)]/\log((\sinh x)/x)$ is strictly increasing on $(0,\pi/2).$  He used this result to prove that the best positive constants $p$ and $q$ for which
$$
\left(\frac{\sinh x}{x}\right)^p<\frac{x}{\sin x}<\left(\frac{\sinh x}{x}\right)^q
$$
is true are $p=1$ and $q=[\log(\pi/2)]/\log(\sinh(\pi/2)/(\pi/2))\approx1.18.$

In an unpublished manuscript, C. Barbu and  L.-I.  Pi\c{s}coran \cite{bapi} have proved, in particular, that
$$
(1-x^2/3)^{-1/4} < \frac{\sinh x}{x}<1+\frac{x^2}{5},\ \ x\in (0,1).
$$

M.-K. Kuo \cite{ko} has developed a method of obtaining an increasing sequence of lower bounds and a decreasing sequence of upper bounds for $(\sin x)/x$, and he has conjectured that the two sequences converge uniformly to $(\sin x)/x$.

Since there is a close connection between the function $(\sin x)/x$ and the Bessel function $J_{1/2}(x)$ (cf. \cite{z10}), it is natural for authors to seek analogs of the Jordan inequality for Bessel and closely related functions.   \'A. Baricz and S.-H. Wu \cite{bar1,bw},  L. Zhu \cite{z4,z10}, and D.-W. Niu et al \cite{ncq} have produced inequalities of this type.  L. Zhu \cite{z1} has also obtained Jordan-type inequalities for $((\sin x)/x)^p$ for any $p>0$.  S.-H. Wu and L. Debnath \cite{wd} have generalized Jordan's inequality to  functions  $f(x)/x$ on $[0,\theta]$ such that  $f$ is $(n+1)$-times differentiable, $f(0)=0,$ and either  $n$ is a positive even integer with  $f^{(n+1)}$  increasing on $[0,\theta]$ or $n$ is a positive odd integer with  $f^{(n+1)}$ decreasing on $[0,\theta].$

\section{Other inequalities involving circular and hyperbolic functions}

\subsection{Redheffer}

In 1968 R. Redheffer \cite{r} proposed the problem of showing that
\begin{equation}\label{Red1}
\frac{\sin \pi x}{\pi x} \ge \frac{1-x^2}{1+x^2},  \ {\rm for\ all\ real}\ x,
\end{equation}
 or, equivalently, that
\begin{equation}\label{Red2}
\frac{\sin x}{x} \ge \frac{\pi^2-x^2}{\pi^2+x^2},  \ {\rm for\ all\ real}\ x.
\end{equation}

A solution of this problem was provided by J. T. Williams \cite{ws}, using infinite products, who also proved the stronger inequality
$$
\frac{\sin \pi x}{\pi x} \ge \frac{1-x^2}{1+x^2} + \frac{(1-x)^2}{x(1+x^2)},  \ \ {\rm for}\  x \ge 1.
$$
Later, using Erd\"os-Tur\'an series and harmonic analysis, J.-L. Li and Y.-L. Li \cite{ll} proved the double inequality
$$
\frac{(1-x^2)(4-x^2)(9-x^2)}{x^6-2x^4+13x^2+36}\le \frac{\sin \pi x}{\pi x}\le \frac{1-x^2}{\sqrt{1+3x^4}},\ \   {\rm for}\ \ 0<x<1.
$$
They also found a method for obtaining new bounds from old for $(\sin x)/x$, but M.-K. Kuo \cite{ko} gave an example to show that the new bounds are not necessarily stronger.
In 2003 C.-P. Chen et al \cite{czq}, using mathematical induction and infinite products, found analogs of the Redheffer inequality for $\cos x$:
$$
\cos x\ge \frac{\pi^2-4x^2}{\pi^2+4x^2},\ \ {\rm for}\ |x|\le \frac{\pi}{2},
$$
and for hyperbolic functions:
$$
\frac{\sinh x}{x}\le\frac{\pi^2+x^2}{\pi^2-x^2},\ \ {\rm for}\ 0 < |x| \le \pi;\ \  {\cosh x}\le\frac{\pi^2+4x^2}{\pi^2-4x^2},\ \ {\rm for}\ |x|\le \frac{\pi}{2}.
$$

In 2008, inspired by the inequalities above, L. Zhu and J.-J. Sun \cite{zs}  proved that
$$
\left( \frac{\pi^2-4x^2}{\pi^2+4x^2}\right)^{\alpha}\le \cos x\le\left( \frac{\pi^2-4x^2}{\pi^2+4x^2}\right)^{\beta}, \ \ \ {\rm for}\ 0\le x \le \frac{\pi}{2},
$$
with best possible constants $\alpha=1$ and $\beta=\pi^2/16$,
and
$$
\left( \frac{\pi^2-x^2}{\pi^2+x^2}\right)^{\gamma}\le \frac{\sin x}{x}\le\left( \frac{\pi^2-x^2}{\pi^2+x^2}\right)^{\delta}, \ \ {\rm for}\ 0<x<\pi,
$$
with best possible constants $\gamma=1$ and $\delta=\pi^2/12$.  They obtained similar results for the hyperbolic sine and cosine functions.
In 2009 L. Zhu \cite{z15} showed that
\begin{equation*}
\left(\frac{\pi^2-x^2}{\sqrt{\pi^4+3x^4}}\right)^{\alpha}\le\frac{\sin x}{x}\le \left(\frac{\pi^2-x^2}{\sqrt{\pi^4+3x^4}}\right)^{\beta}, \ \ 0<x\le\pi,
\end{equation*}
holds if and only if $\alpha \ge\pi^2/6$ and $\beta\le 1$, with analogous results for $\cos x$ and $(\tan x)/x$.  In 2009 \'A. Baricz and S.-H. Wu \cite{bw2} and in 2011 L. Zhu \cite{z16}  proved Redheffer-type inequalities for Bessel functions.

\subsection{Cusa-Huygens}

The inequality
$$
\frac{\sin x}{x}<\frac{2+\cos x}{3},\ \ 0<x<\pi/2,
$$
was discovered by N. de Cusa in the fifteenth century (cf. \cite{cc}), and proved rigorously by Huygens \cite{hs} in the seventeenth century.  In 2009 L. Zhu
\cite{z12} obtained the following inequalities of Cusa-Huygens type:
$$
\left(\frac{\sin x}{x}\right)^{\alpha}<\frac{2}{3}+\frac{1}{3}(\cos x)^{\alpha}, \ \ 0<x<\frac{\pi}{2},\ \ \alpha\ge 1,
$$
and
$$\left(\frac{\sinh x}{x}\right)^{\alpha}<\frac{2}{3}+\frac{1}{3}(\cosh x)^{\alpha},\ \ x>0,\ \ \alpha \ge 1.
$$

That same year L. Zhu \cite{z6} discovered a more general set of inequalities of Cusa type, from which many other types of inequalities for circular functions can be derived.  He proved the following:  Let  $0<x<\pi/2$.  If $p\ge1$, then
\begin{equation}\label{source}
(1-\alpha )+\alpha(\cos x)^p<\left(\frac{\sin x}{x}\right)^p<(1-\beta)+\beta (\cos x)^p
\end{equation}
if and only if  $\beta \le 1/3$ and $\alpha\ge 1-(2/\pi)^p.$  If  $0\le p\le 4/5$, then \eqref{source} holds if and only if $\alpha\ge 1/3$ and $\beta\le 1-(2/\pi)^p.$  If $p<0$, then the second inequality in \eqref{source} holds if and only if $\beta\ge 1/3.$  In a later paper \cite{z10} L. Zhu  obtained estimates for $(\sin x)/x$ and $(\sinh x)/x$ that led to new infinite series for these functions. For some similar results see also \cite{wb}.

In 2011 C.-P. Chen and W.-S. Cheung \cite{cc} obtained the sharp Cusa-Huygens-type inequality
\begin{equation*}
\left(\frac{2+\cos x}{3}\right)^{\alpha}<\frac{\sin x}{x}<\left(\frac{2+\cos x}{3}\right)^{\beta},
\end{equation*}
for $0<x<\pi/2$, with best possible constants  $\alpha = (\log(\pi/2))/\log(3/2)\approx 1.11$ and $\beta=1.$

In 2011 E. Neuman and J. S\'andor \cite{ns} discovered a pair of optimal inequalities for hyperbolic and trigonometric functions, proving that, for $0<x<\pi/2$, the best positive constants $p$ and $q$ in the inequality
\begin{equation*}
\frac{1}{(\cosh x)^p}<\frac{\sin x}{x}<\frac{1}{(\cosh x)^q}
\end{equation*}
are $p=(\log(\pi/2))/\log\cosh (\pi/2)\approx 0.49$ and $q=1/3$,
 and that for $x\ne 0$ the best positive constants $p$ and $q$ in the inequality
\begin{equation*}
\left(\frac{\sinh x}{x}\right)^p<\frac{2}{\cos x + 1}< \left(\frac{\sinh x}{x}\right)^q
\end{equation*}
are $p=3/2$ and $q=(\log 2)/\log[(\sinh (\pi/2))/(\pi/2)]\approx 1.82.$

\subsection{Becker-Stark}

In 1978 M. Becker and E. L. Stark \cite{bes} obtained the double inequality
$$
\frac{8}{\pi^2-4x^2}<\frac{\tan x}{x}<\frac{\pi^2}{\pi^2-4x^2},\ \  0<x<\frac{\pi}{2},
$$
where the numerator constants $8$ and $\pi^2$ are best possible.

In 2008 L. Zhu and J.-J. Sun \cite{zs}  showed that
\begin{equation*}
\left(\frac{\pi^2+4x^2}{\pi^2-4x^2}\right)^{\alpha}\le\frac{\tan x}{x}\le\left(\frac{\pi^2+4x^2}{\pi^2-4x^2}\right)^{\beta}, \ \ 0<x<\frac{\pi}{2},
\end{equation*}
holds if and only if $\alpha\le\pi^2/24$ and $\beta\ge 1$.

In 2010 L. Zhu and J.-K. Hua \cite{zh} sharpened the Becker-Stark inequality by proving that
\begin{equation*}
\frac{\pi^2+\alpha x^2}{\pi^2-4x^2}<\frac{\tan x}{x}<\frac{\pi^2+\beta x^2}{\pi^2-4x^2}, \ \ 0<x<\frac{\pi}{2},
\end{equation*}
where $\alpha=4(8-\pi^2)/\pi^2\approx -0.76$ and $\beta=\pi^2/3-4\approx-0.71$ are the best possible constants.  They also developed a systematic method for obtaining a sequence of sharp inequalities of this sort.

In 2011  H.-F. Ge \cite{g} obtained
$$
\frac{8}{\pi^2-4x^2}+\left(1-\frac{8}{\pi^2}\right)<\frac{\tan x}{x}<\frac{\pi^4}{12}\frac{1}{\pi^2-4x^2}+\left(1-\frac{\pi^2}{12}\right),
$$
for $0<x<\pi/2$.  That same year C.-P. Chen and W.-S. Cheung \cite{cc} proved the sharp Becker-Stark-type inequality
\begin{equation*}
\left(\frac{\pi^2}{\pi^2-4x^2}\right)^{\alpha}<\frac{\tan x}{x}< \left(\frac{\pi^2}{\pi^2-4x^2}\right)^{\beta},
\end{equation*}
with best possible constants $\alpha = \pi^2/12\approx 0.82$ and $\beta=1.$

\subsection{Wilker}

In 1989 J. Wilker \cite{wr} posed the problem of proving that
\begin{equation}\label{Wilker1}
\left(\frac{\sin x}{x}\right)^2 +\frac{\tan x}{x}>2, \ \ {\rm for}\ 0<x<\frac{\pi}{2},
\end{equation}
and of finding
\begin{equation}\label{WilkerInf}
c\equiv \inf_{0<x<\pi/2}\frac{\left(\frac{\sin x}{x}\right)^2+\frac{\tan x}{x}-2}{x^3\tan x}.
\end{equation}
J. Anglesio et al \cite{sjva} showed that the function in \eqref{WilkerInf}  is decreasing on $(0,\pi/2)$, that the value of $c$ in \eqref{WilkerInf}  is $16/\pi^4$, and that, moreover, the supremum of the expression in  \eqref{WilkerInf}  on $(0,\pi/2)$ is $8/45$.   Hence
\begin{equation}\label{Wilker2}
2+\frac{16}{\pi^4}x^3\tan x\le\left(\frac{\sin x}{x}\right)^2+\frac{\tan x}{x}\le2+\frac{8}{45}x^3\tan x,
\end{equation}
for $0 < x < \pi/2$, where  $16/\pi^4\approx 0.164$ and $8/45\approx 0.178$ are best possible constants.  (Note:   \cite{avv4} erroneously quoted \cite{sjva} as saying that the function in \eqref{WilkerInf}  is increasing.)
In 2007 S.-H. Wu and H. M. Srivastava \cite{ws1} proved the Wilker-type inequality
\begin{equation}\label{BS}
\left(\frac{x}{\sin x}\right)^2 +\frac{x}{\tan x}>2, \ \ {\rm for}\ 0<x<\frac{\pi}{2}.
\end{equation}
However, \'A. Baricz and J. S\'andor \cite{bs} discovered that \eqref{BS}  is implied by \eqref{Wilker1}.

In 2009 L. Zhu \cite{z12} generalized \eqref{Wilker1} and obtained analogs for hyperbolic functions, showing that,  for $0<x<\pi/2$, $\alpha\ge 1,$
\begin{equation*}
\left(\frac{\sin x}{x}\right)^{2\alpha}+\left(\frac{\tan x}{x}\right)^{\alpha}> \left(\frac{x}{\sin x}\right)^{2\alpha}+\left(\frac{x}{\tan x}\right)^{\alpha}>2,
\end{equation*}
and that, for $x>0$, $\alpha\ge 1$,
\begin{equation*}
\left(\frac{\sinh x}{x}\right)^{2\alpha}+\left(\frac{\tanh x}{x}\right)^{\alpha}> \left(\frac{x}{\sinh x}\right)^{2\alpha}+\left(\frac{x}{\tanh x}\right)^{\alpha}>2.
\end{equation*}
These two results of Zhu are special cases of a recent lemma due to E. Neuman \cite[Lemma 2]{n4}.

In 2012 J. S\'andor \cite{s6} has proved that, for $0<x\leq\pi/2$, $\alpha> 0,$
$$
\left(\frac{x}{\sin x}\right)^{2\alpha}+\left(\frac{x}{\sinh x}\right)^{\alpha}>
\left(\frac{\sinh x}{x}\right)^{2\alpha}+\left(\frac{\sin x}{x}\right)^{\alpha}>2.
$$

Using power series, C.-P. Chen and W.-S. Cheung \cite{cc2} obtained the following sharper versions of \eqref{Wilker2}:
\begin{equation}\label{cc1}
\frac{16}{315}x^5\tan x<\left(\frac{\sin x}{x}\right)^2+\frac{\tan x}{x}-\left[2+\frac{8}{45}x^4\right]<\left(\frac{2}{\pi}\right)^6x^5\tan x,
\end{equation}
and
\begin{equation}\label{cc2}
\frac{104}{4725}x^7\tan x<\left(\frac{\sin x}{x}\right)^2+\frac{\tan x}{x}-\left[2+\frac{8}{45}x^4+\frac{16}{315}x^6\right]<\left(\frac{2}{\pi}\right)^8x^7\tan x.
\end{equation}
The constants $16/315\approx 0.051$ and $(2/\pi)^6\approx 0.067$ in \eqref{cc1} and $104/4725\approx 0.022$ and $(2/\pi)^8\approx 0.027$ in \eqref{cc2} are best possible.
For $0<x<\pi/2$,   Chen and Cheung also obtained upper estimates complementary to \eqref{BS}:
\begin{equation*}
\left(\frac{x}{\sin x}\right)^2 +\frac{x}{\tan x}<2+\frac{2}{45}x^3\tan x
\end{equation*}
and
\begin{equation*}
\left(\frac{x}{\sin x}\right)^2 +\frac{x}{\tan x}<2+\frac{2}{45}x^4+\frac{8}{945}x^5\tan x,
\end{equation*}
where the constants  $2/45$ and $8/945$ are best possible.

In 2012, J. S\'andor \cite{s1} has shown that
\begin{equation*}
\frac{\sin x}{x}+q\frac{\sinh x}{x}>q+1,\ x\ne 0,
\end{equation*}
and
\begin{equation*}
\left(\frac{\sinh x}{x}\right)^q+\frac{\sin x}{x}>2, \ 0<x<\frac{\pi}{2},
\end{equation*}
where $q=[\log(\pi/2)]/\log[(\sinh(\pi/2))/(\pi/2)]\approx 1.18.$

Extensions of the generalized Wilker inequality for Bessel functions were obtained by \'A. Baricz and J. S\'andor \cite{bs} in 2008.

\subsection{Huygens}

An older inequality  due to C. Huygens \cite{hs} is similar in form to \eqref{Wilker1}:
\begin{equation}\label{Huy}
2\left(\frac{\sin x}{x}\right)+\frac{\tan x}{x} >3, \ \ {\rm for}\ 0<|x|<\frac{\pi}{2},
\end{equation}
and actually implies \eqref{Wilker1} (see \cite{ns1}).  In 2009 L. Zhu \cite{z5} obtained the following inequalities of Huygens type:
\begin{equation*}
(1-p)\frac{\sin x}{x}+p\frac{\tan x}{x}>1>(1-q)\frac{\sin x}{x}+q\frac{\tan x}{x}
\end{equation*}
for all $x\in(0,\pi/2 )$ if and only if $p\ge 1/3$ and $q\le 0$;
\begin{equation*}
(1-p)\frac{\sinh x}{x}+p\frac{\tanh x}{x}>1>(1-q)\frac{\sinh x}{x}+q\frac{\tanh x}{x}
\end{equation*}
for all $x\in (0,\infty )$ if and only if $p\le 1/3$ and $q\ge 1$;
\begin{equation*}
(1-p)\frac{x}{\sin x}+p\frac{x}{\tan x}>1>(1-q)\frac{x}{\sin x}+q\frac{x}{\tan x}
\end{equation*}
for all $x\in (0,\pi/2)$ if and only if $p\le 1/3$ and $q\ge1-2/\pi$; and
$$
(1-p)\frac{x}{\sinh x} + p\frac{x}{\tanh x} > 1 > (1-q)\frac{x}{\sinh x}+q\frac{x}{\tanh x}
$$
for all $x\in (0,\infty )$ if and only if $p\ge 1/3$ and $q\le 0$.

In 2012 J. S\'andor \cite{s6} has showed that, for $0<x\leq\pi/2$, $\alpha>0$,
$$
2\left(\frac{\sinh x}{x}\right)^{\alpha}+\left(\frac{\sin x}{x}\right)^{\alpha}>
2\left(\frac{x}{\sin x}\right)^{\alpha}+\left(\frac{x}{\sinh x}\right)^{\alpha}>3.
$$

C.-P. Chen and W.-S. Cheung \cite{cc2} also found sharper versions of \eqref{Huy}, as follows:
For $0<x<\pi/2$,
\begin{equation}\label{Huy1}
3+\frac{3}{20}x^3\tan x<2\left(\frac{\sin x}{x}\right)+\frac{\tan x}{x}<3+\left(\frac{2}{\pi}\right)^4x^3\tan x
\end{equation}
and
\begin{equation}\label{Huy2}
\frac{3}{56}x^5\tan x<2\left(\frac{\sin x}{x}\right)+\frac{\tan x}{x}-\left[3+\frac{3}{20}x^4\right]<\left(\frac{2}{\pi}\right)^6x^5\tan x,
\end{equation}
where the constants $3/20=0.15$ and $(2/\pi)^4\approx 0.16$ in \eqref{Huy1} and $3/56\approx 0.054$ and $(2/\pi )^6\approx 0.067$ in \eqref{Huy2} are best possible.

Recently Y. Hua \cite{hu} have proved the following sharp inequalities:
For $0<|x|<\pi/2$,
$$
3+\frac{1}{40}x^3\sin x<\frac{\sin x}{x}+2\frac{\tan(x/2)}{x/2}<3+\frac{80-24\pi}{\pi^4}x^3\sin x,
$$
where the constants $1/40$ and $(80-24\pi)/\pi^4$ are best possible;
and, for $x\neq0$,
$$
3+\frac{3}{20}x^3\tanh x<2\frac{\sinh x}{x}+\frac{\tanh x}{x}<3+\frac{3}{20}x^3\sinh x,
$$
where the constant $3/20$ is best possible.

\subsection{Shafer}

The problem of proving
$$
\arctan x> \frac{3x}{1+2\sqrt{1+x^2}}, \ \ x>0,
$$
was proposed by R. E. Shafer in 1966.  Solutions were obtained by L. S. Grinstein, D. C. B. Marsh, and J. D. E. Konhauser \cite{sgmk} in 1967.  In 2011 C.-P. Chen, W.-S. Cheung, and W.-S. Wang \cite{ccw} found, for each $a>0$, the largest number $b$ and the smallest number $c$ such that the inequalities
$$
\frac{bx}{1+a\sqrt{1+x^2}}\le \arctan x \le \frac{cx}{1+a\sqrt{1+x^2}}
$$
are valid for all $x\ge 0$.  Their answer to this question is indicated in the following table:

\begin{center}
\begin{tabular}{|c| c |c|}
\hline
$a$&\rm{largest} \ $ b$&\rm{smallest} \ $c$\\ \hline
$0<a\le\pi/2$&$b=\pi a/2$&$c=1+a$\\ \hline
$\pi/2<a\le2/(\pi-2)$&$b=4(a^2-1)/a^2$&$c=1+a$\\ \hline
 $2/(\pi-2)<a<2$&$ b=4(a^2-1)/a^2$&$c=\pi a/2$\\ \hline
$2\le a<\infty$&$ b=1+a$&$c=\pi a/2$\\
\hline
\end{tabular}
\end{center}

\medskip

In 1974, in a numerical analytical context \cite{sh2},  R. E. Shafer presented the inequality
\begin{equation*}
\arctan x\ge \frac{8x}{3+\sqrt{25+(80/3)x^2}},\ x> 0,
\end{equation*}
which he later proved analytically \cite{sh3}.
In \cite{z7} L. Zhu proved that the constant $80/3$ in Shafer's inequality is best possible, and also obtained the complementary inequality
\begin{equation*}
\arctan x < \frac{8x}{3+\sqrt{25+(256/\pi^2)x^2}},\ x>0,
\end{equation*}
where $256/\pi^2$ is the best possible constant.

\subsection{Fink}

In \cite[p. 247]{m}, there is a lower bound for $\arcsin x$ on $[0,1]$ that is similar to Shafer's for $\arctan x$.  In 1995 A. M. Fink \cite{f} supplied a complementary upper bound.  The resulting double inequality is
\begin{equation}\label{eq:fink}
\frac{3x}{2+\sqrt{1-x^2}}\le \arcsin x\le\frac{\pi x}{2+\sqrt{1-x^2}}, \ 0\le x\le1,
\end{equation}
and both numerator constants are best possible.  Further refinements of these inequalities, along with analogous ones for ${\rm arcsinh}\,x$ were obtained by  L. Zhu \cite{z13} and by W.-H. Pan with L. Zhu \cite{pz}.  We note that, for $0<x<1$, this double inequality  is equivalent to
$$
\frac{2+\cos x}{\pi} < \frac{\sin x}{x} < \frac{2+\cos x}{3},\ \ 0<x<\frac{\pi}{2},
$$
in which the second relation is the Cusa inequality.

\subsection{Carlson}

In 1970 B. C. Carlson \cite[(1.14)]{c} proved the inequality
\begin{equation}\label{carlson}
\frac{6\sqrt{1-x}}{2\sqrt{2}+\sqrt{1+x}}<\arccos x< \frac{\sqrt[3]{4}\cdot\sqrt{1-x}}{(1+x)^{1/6}},\ \ 0\le x < 1.
\end{equation}
In 2012, seeking to sharpen and generalize \eqref{carlson}, C.-P. Chen and C. Mortici \cite{cm} determined, for each fixed $c>0$,  the largest number $a$ and smallest number $b$ such that the double inequality
$$
\frac{a\sqrt{1-x}}{c+\sqrt{1+x}}\le\arccos x\le \frac{b\sqrt{1-x}}{c+\sqrt{1+x}}
$$
is valid for all $x\in [0,1]$.
Their answer to this question is indicated in the following table:

\medskip

\begin{center}
\begin{tabular}{|c| c |c|}
\hline
$c$&\rm{largest} \ $ a$&\rm{smallest} \ $b$\\ \hline
$0<x<(2\pi-4)/(4-\pi)$&$(1+a)\pi/2$&$2+\sqrt{2}a$\\ \hline
$(2\pi-4)/(4-\pi)\leq x\leq(4-\pi)/(\pi-2\sqrt{2})$&$8(a^2-2)/a^2$&$2+\sqrt{2}a$\\ \hline
$(4-\pi)/(\pi-2\sqrt{2})<x<2\sqrt{2}$&$4(a^2-1)/a^2$&$(1+a)\pi/2$\\ \hline
$2\sqrt{2}\leq x<\infty$&$2+\sqrt{2}a$&$(1+a)\pi/2$\\
\hline
\end{tabular}
\end{center}

\medskip

\noindent These authors also proved that, for all $x\in [0,1]$, the inequalities
$$
\frac{\sqrt[3]{4}\cdot\sqrt{1-x}}{a+(1+x)^{1/6}}\le \arccos x \le \frac{\sqrt[3]{4}\cdot\sqrt{1-x}}{b+(1+x)^{1/6}}
$$
hold on $[0,1]$, with best constants $a=(2\sqrt[3]{4}-\pi)/\pi\approx 0.01$ and $b=0.$

Moreover, in view of the right side of \eqref{carlson}, in 2011 C.-P. Chen, W.-S. Cheung, and W.-S. Wang \cite{ccw}
considered functions of the form
$$
f(x)\equiv \frac{r(1-x)^p}{(1+x)^q}
$$
on $[0,1]$, and determined the values of $p, q, r$ such that $f(x)$ is the best $3$rd-order approximation of $\arccos x$ in a neighborhood of the origin.  The answer is that, for $p=(\pi+2)/\pi^2$, $q=(\pi-2)/\pi^2$, $r = \pi/2$, one has
\begin{equation*}
\lim_{x\to 0}\frac{\arccos x-f(x)}{x^3}=\frac{\pi^2-8}{6\pi^2}\ .
\end{equation*}
With the values of $p, q, r$ stated above, the authors were led to a new lower bound for $\arccos $:
\begin{equation*}
\arccos x\ge\frac{(\pi/2)(1-x)^{(\pi+2)/\pi^2}}{(1+x)^{(\pi-2)/\pi^2}},\ \ \ 0<x\le1.
\end{equation*}

\subsection{Lazarevi\'c}

In \cite{l} I. Lazarevi\'c proved that, for $x\ne 0,$
$$
\left(\frac{\sinh x}{x}\right)^q> \cosh x
$$
if and only if $q\ge 3.$   L. Zhu improved upon this inequality in \cite{z14}, by showing that if $p>1$ or $p\le8/15$ then
$$
\left(\frac{\sinh x}{x}\right)^q>p+(1-p)\cosh x
$$
for all  $x>0$  if and only if $q\ge3(1-p).$
For some similar results see also \cite{wb}.

In 2008  \'A. Baricz \cite{bar6} extended the Lazarevi\'c inequality to modified Bessel functions and also deduced some Tur\'an- and Lazarevi\'c-type inequalities for the confluent hypergeometric functions.

\subsection{Neuman}

E. Neuman \cite{n} has recently established several inequalities involving new combinations of circular and hyperbolic functions.  In particular, he has proved that if $x\ne 0,$ then
\begin{equation*}
(\cosh x)^{2/3}<\frac{\sinh x}{\arcsin (\tanh x)}<\frac{1+2\cosh x}{3},
\end{equation*}
\begin{equation*}
[(\cosh 2x)^{1/2}\cosh^2 x]^{1/3}<\frac{\sinh x}{\rm{arcsinh} (\tanh x)}<\frac{(\cosh 2x)^{1/2}+2\cosh x}{3},
\end{equation*}
and
\begin{equation*}
[(\cosh 2x)\cosh x]^{1/3}<\frac{\sinh x}{\arctan (\tanh x)}<\frac{2(\cosh 2x)^{1/2}+\cosh x}{3}.
\end{equation*}

\section{Euler's gamma function}

For $\Re z>0$ the \emph{gamma function} is defined by
$$
\Gamma(z)\equiv\int_0^{\infty}t^{z-1}e^{-t}dt,
$$
and the definition  is extended by analytic continuation to the entire complex plane minus the set of nonpositive integers.  This function, discovered by Leonhard Euler in 1729, is a natural generalization of the factorial, because of the functional identity
$$
\Gamma(z+1)=z\Gamma(z).
$$
The gamma function is one of the best-known and most important special functions in mathematics, and has been studied intensively.

We begin our treatment of this subject by considering an important special constant discovered by Euler and related to the gamma function.

\subsection{The Euler-Mascheroni constant and harmonic numbers}

The \emph{Euler-Mascheroni constant} $\gamma=0.5772156649...$ is defined as
\begin{equation}\label{eulermascheroni}
\gamma\equiv\lim_{n\to\infty}\gamma_n,
\end{equation}
where $\gamma_n \equiv H_n-\log n$, $n\in\mathbb{N},$ and where $H_n$ are the \emph{harmonic numbers}
\begin{equation}\label{harmonic}
H_n\equiv \sum_{k=1}^n\frac{1}{k}=\int_0^1\frac{1-x^n}{1-x}dx.
\end{equation}
The number $\gamma$ is one of the most important constants in mathematics, and is useful in analysis, probability theory, number theory, and other branches of pure and applied mathematics.  The numerical value of $\gamma$ is known to $29,844,489,545$ decimal places, thanks to computation by Yee and Chan in 2009 \cite{ye} (see \cite[p. 273]{chl}).

The sequence $\gamma_n$ converges very slowly to $\gamma,$ namely with order $1/n.$  By replacing $\log n$ in this sequence by $\log (n+1/2)$, D. W. DeTemple \cite{d} obtained quadratic convergence (see also \cite{cn2}).   In \cite{m4} C. Mortici made a careful study of how convergence is affected by changes in the logarithm term.  He introduced new sequences
$$
M_n\equiv H_n-\log\frac{P(n)}{Q(n)},
$$
where  $P$ and $Q$ are polynomials with leading coefficient $1$  and $\deg P - \deg Q = 1.$   By judicious choice of the degrees and coefficients of $P$ and $Q$ he was able to produce sequences $M_n$ tending to $\gamma$ with convergence of order $1/n^4$ and $1/n^6$.  He also gave a recipe for obtaining sequences converging to $\gamma$ with order $1/n^{2k+2},$ where $k$ is any positive integer.  This study is based on the author's lemma, proved in \cite{m5}, that connects the rate of convergence of a convergent sequence $\{x_n\}$ to that of the sequence $\{x_n-x_{n+1}\}$.

In 1997  T. Negoi \cite{ni} showed that if $T_n\equiv H_n-\log(n+1/2+1/(24n))$, then $T_n + [4n^3]^{-1}$ is strictly decreasing to $\gamma$ and $T_n + [48(n+1)]^{-3}$ is strictly increasing to $\gamma$, so that $[48(n+1)]^{-3}<\gamma-T_n<[48n^3]^{-1}. $  In 2011 C.-P. Chen \cite{cn} established sharper bounds for $\gamma-T_n$, by using a lemma of Mortici \cite{m5}.

Using another approach, in 2011 E. Chlebus \cite{chl} developed a recursive scheme for modifying the sequence $H_n-\log n$ to accelerate the convergence to $\gamma$  to any desired order.  The first step in Chlebus' scheme is equivalent to the DeTemple \cite{d} approximation, while the next step yields a sequence that closely resembles the one due to T. Negoi \cite{ni}.

In \cite{a} H. Alzer studied the harmonic numbers \eqref{harmonic}, obtaining several new inequalities for them.  In particular, for $n\ge 2$, he proved that
\begin{equation*}
\alpha\frac{\log (\log n+\gamma )}{n^2}\le H_n^{1/n}-H_{n+1}^{1/(n+1)}<\beta\frac{\log (\log n+\gamma )}{n^2},
\end{equation*}
where $\alpha=(6\sqrt{6}-2\sqrt[3]{396})/(3\log (\log 2+\gamma ))\approx  0.014$ and $\beta = 1$ are the best possible constants and $\gamma$ is the Euler-Mascheroni constant.

\subsection{Estimates for the gamma function}

In \cite[Lemma 2.39]{avv6} G. D. Anderson, M. Vamanamurthy, and M. Vuorinen proved that
\begin{equation}\label{AQ}
\lim_{x\to\infty}\frac{\log \Gamma\left(\frac{x}{2}+1\right)}{x\log x}=\frac{1}{2}
\end{equation}
and that the function $(\log \Gamma (1+x/2))/x$ is strictly increasing from $[2,\infty)$ onto $[0,\infty ).$
In \cite{aq} G. D. Anderson and S.-L. Qiu showed that $(\log\Gamma (x+1))/(x\log x)$  is strictly increasing from $(1,\infty)$ onto
$(1-\gamma,1),$ where $\gamma$ is the Euler-Mascheroni constant defined by \eqref{eulermascheroni}, thereby obtaining the strict inequalities
\begin{equation}\label{gammaAQ}
x^{(1-\gamma)x-1}<\Gamma (x)<x^{x-1},\ \ x>1 .
\end{equation}
They also conjectured that the function $(\log \Gamma(x+1))/(x\log x)$ is concave on $(1,\infty)$, and this conjecture was proved by \'A. Elbert and A. Laforgia in \cite[Section 3]{el}. One should note that in 1989 J. S\'andor \cite{s7}
proved that the function $(\Gamma(x+1))^{1/x}$ is strictly concave for  $x\geq7$.

Later H. Alzer \cite{a1} was able to extend \eqref{gammaAQ} by proving that, for $x\in(0,1)$,
\begin{equation}\label{gammaAlzer}
x^{\alpha (x-1)-\gamma}<\Gamma (x)<x^{\beta (x-1)-\gamma},
\end{equation}
with best possible constants $\alpha=1-\gamma= 0.42278\dots$ and $\beta=(\pi^2/6-\gamma)/2=0.53385\ldots$.
For $x\in (1,\infty)$ H. Alzer was able to sharpen \eqref{gammaAQ} by showing that \eqref{gammaAlzer} holds with best possible constants $\alpha=(\pi^2/6-\gamma)/2\approx 0.534$ and $\beta=1$.  His principal new tool was the convolution theorem for Laplace transforms.

Another type of approximation for $\Gamma(x)$ was derived by P. Iv\'ady \cite{i} in 2009:
\begin{equation}\label{gammaIvady}
\frac{x^2+1}{x+1}<\Gamma(x+1)<\frac{x^2+2}{x+2},\ \ 0 < x < 1.
\end{equation}
In 2011 J.-L. Zhao, B.-N. Guo, and F. Qi \cite{zgq} simplified and sharpened \eqref{gammaIvady} by proving that the function
$$
Q(x)\equiv \frac{\log\Gamma(x+1)}{\log(x^2+1)-\log(x+1)}
$$
is strictly increasing from $(0,1)$ onto $(\gamma,2(1-\gamma))$, where $\gamma$ is the Euler-Mascheroni constant.  As a consequence, they proved that
\begin{equation*}
\left(\frac{x^2+1}{x+1}\right)^{\alpha}< \Gamma(x+1)<\left(\frac{x^2+1}{x+1}\right)^{\beta}, \ \ 0 < x < 1,
\end{equation*}
with best possible constants $\alpha=2(1-\gamma)$ and $\beta=\gamma.$

Very recently C. Mortici \cite{m8} has determined by numerical experiments that the upper estimate in (\ref{gammaAQ}) is a better approximation for  $\Gamma (x)$  than the lower one when  $x$ is very large.  Hence, he has sought estimates of the form $\Gamma (x)\approx x^{a(x)}$, where $a(x)$ is close to $x-1$ as $x$ approaches infinity.  For example, he proves that
$$
x^{(x-1)a(x)}<\Gamma (x) < x^{(x-1)b(x)}, \quad x\ge 20,
$$
where $a(x)=1-1/\log x+1/(2x)-(1-(\log 2\pi)/2)/(x\log x)$ and where $b(x)=1-1/\log x+1/(2x)$.  The left inequality is valid for $x\ge 2$.  Mortici has also obtained a pair of sharper inequalities of this type, valid for $x\ge 2$, and has showed how lower and upper estimates of any desired accuracy may be obtained.
His proofs are based on an approximation for $\log\Gamma (x)$  in terms of series involving Bernoulli numbers \cite[p. 29]{aar} and on truncations of an asymptotic series for the function
$(\log \Gamma (x))/((x-1)\log x).$  These results provide improvements of (\ref{gammaAQ}).

\subsection{Factorials and Stirling's formula}

The well-known  \emph{Stirling's formula} for $n!$,
\begin{equation}\label{stirling}
\alpha_n\equiv \left(\frac{n}{e}\right)^n\sqrt{2\pi n}\, ,
\end{equation}
discovered by the precocious home-schooled and largely self-taught eighteenth century Scottish mathematician James Stirling, approximates $n!$ asymptotically in the sense that
$$
\lim_{n\to\infty}\frac{n!}{\alpha_n}=1.
$$
 Because of the importance of this formula in probability and statistics, number theory, and scientific computations, several authors have sought to replace \eqref{stirling} by a simple sequence that approximates $n!$ more closely (see the discussions in \cite{bat} and \cite{bat4}).  For example, W. Burnside \cite{bu} proved in 1917 that
\begin{equation}\label{burnside}
n!\sim  \beta_n \equiv \sqrt{2\pi}\left(\frac{n+1/2}{e}\right)^{n+1/2},
\end{equation}
that is, $\lim\limits_{n\to\infty}(n!/\beta_n)=1.$ In 2008 N. Batir \cite{bat} determined that the best constants  $a$ and $b$ such that
\begin{equation}\label{batir}
\frac{n^{n+1}e^{-n}\sqrt{2\pi}}{\sqrt{n-a}}\le n!< \frac{n^{n+1}e^{-n}\sqrt{2\pi}}{\sqrt{n-b}}
\end{equation}
are $a=1-2\pi e^{-2}\approx 0.1497$ and $b=1/6\approx 0.1667$.  Batir offers a numerical table illustrating that his upper bound formula $n^{n+1}e^{-n}\sqrt{2\pi}/\sqrt{n-1/6}$ gives much better approximations to $n!$ than does either \eqref{stirling} or \eqref{burnside}.

In a later paper \cite{bat4} N. Batir observed that many of the improvements of Stirling's formula take the form
\begin{equation}\label{batir2}
n!\sim e^{-a}\left(\frac{n+a}{e}\right)^n\sqrt{2\pi (n+b)}
\end{equation}
for some real numbers $a$ and $b$.  Batir sought the pair of constants $a$ and $b$ that would make \eqref{batir2} optimal.  He proved that the best pairs $(a,b)$ are $(a_1,b_1)$ and $(a_2,b_2)$, where
$$
a_1=\frac{1}{3} +\frac{\lambda}{6}-\frac{1}{6}\sqrt{6-\lambda^2+4/\lambda}\approx 0.54032,\ \ b_1=a_1^2+1/6\approx 0.45861
$$
and
$$
a_2=\frac{1}{3}+\frac{\lambda}{6}+\frac{1}{6}\sqrt{6-\lambda^2+4/\lambda}\approx 0.95011,\ \  b_2=a_2^2+1/6\approx 1.06937,
$$
where $\lambda = \sqrt{2+2^{2/3}+2^{4/3}}\approx 2.47128$ and $a_1$ and $a_2$ are the real roots of the quartic equation $3x^4-4x^3+x^2+1/12=0$.

S. Ramanujan \cite{ra} sought to improve Stirling's formula \eqref{stirling} by replacing $\sqrt{2n}$ in the formula by the sixth root of a cubic polynomial in $n$:
\begin{equation}\label{ramastirling}
\Gamma(n+1) \approx \sqrt{\pi}\left(\frac{n}{e}\right)^n\sqrt[6]{8n^3+4n^2+n+\frac{1}{30}}.
\end{equation}
In this connection there appears in the record also his double inequality, for $x\ge 1$,
\begin{equation}\label{ramaineq}
\sqrt[6]{8x^3+4x^2+x+\frac{1}{100}}<\frac{\Gamma(x+1)}{\sqrt{\pi}\left(\frac{x}{e}\right)^x}<\sqrt[6]{8x^3+4x^2+x+\frac{1}{30}}.
\end{equation}
 Motivated by this inequality of Ramanujan, the authors of \cite{avv2} defined the function $h(x)\equiv u(x)^6-(8x^3+4x^2+x)$, where
$u(x)=(e/x)^x\Gamma (x+1)/\sqrt{\pi},$ and conjectured that $h(x)$ is increasing from $(1,\infty)$ into $(1/100,1/30)$.    In 2001 E. A. Karatsuba \cite{ka}
settled this conjecture by showing that  $h(x)$ is increasing from $[1,\infty)$ onto $[h(1),1/30),$ where $h(1)=e^6/\pi^3-13\approx 0.011$.

In an unpublished document, E. A. Karatsuba suggested modifying Ramanujan's approximation formula (\ref{ramastirling}) by replacing the radical with the $2k^{\rm th}$ root of a polynomial of degree $k$, and determining the best such asymptotic approximation.  Such a program was partially realized by C. Mortici \cite{m} in 2011, who proposed formula (\ref{ramamortici}) below for $k=4$, but the more general problem suggested by Karatsuba remains an open problem.  Mortici's proposed Ramanujan-type asymptotic approximation is as follows:
\begin{equation}\label{ramamortici}
\Gamma(n+1)\approx \sqrt{\pi}\left(\frac{n}{e}\right)^n\sqrt[8]{16n^4+\frac{32}{3}n^3+\frac{32}{9}n^2+\frac{176}{405}n-\frac{128}{1215}}\ .
\end{equation}
In connection with \eqref{ramamortici}, he defined the function
$$
g(x)\equiv u(x)^8-\left(16x^4+\frac{32}{3}x^3+\frac{32}{9}x^2+\frac{176}{405}x\right),
$$
where $u(x)=(e/x)^x\Gamma (x+1)/\sqrt{\pi},$ and proved that $g(x)$ is strictly decreasing from $[3,\infty)$ onto $(g(\infty ),g(3)]$, where  $g(\infty )=-128/1215 \approx -0.105$ and $g(3)=256e^{24}/(43046721\pi^4) - 218336/135\approx -0.088.$
Mortici's method for proving monotonicity was simpler than Karatsuba's, because he employed an excellent result of H. Alzer \cite{alzer} concerning complete monotonicity (see section \ref{completemonot} below for definitions).  Mortici claimed that his method would also simplify Karatsuba's proof in \cite{ka}.
Finally, he proved that, for $x\ge3$,
$$
R(x,\alpha )<\frac{\Gamma(x+1)}{\sqrt{\pi}\left(\frac{x}{e}\right)^x}\le
R(x,\beta ),
$$
where $R(x,t)\equiv \sqrt[8]{16x^4+\frac{32}{3}x^3+\frac{32}{9}x^2+\frac{176}{405}x-t}$, and $\alpha=128/1215$, $\beta=g(3)$ are the best possible constants.

In 2012 M. Mahmoud, M. A. Alghamdi, and R. P. Agarwal \cite{maa} deduced a new family of upper bounds for $\Gamma(n+1)$ of the form
$$
\Gamma(n+1)<\sqrt{2\pi n}\left(\frac{n}{e}\right)^n e^{M_n^{[m]}},\ \ n\in\mathbb{N},
$$
$$
M_n^{[m]}\equiv\frac{1}{2m+3}\left[\frac{1}{4n}+\sum_{k=1}^m\frac{2m-2k+2}{2k+1}2^{-2k}\zeta(2k,n+1/2)\right],\ \ n\in\mathbb{N},
$$
where $\zeta$ is the \emph{Hurwitz zeta function}
$$
\zeta(s,q)\equiv\sum_{k=0}^\infty\frac{1}{(k+q)^s}.
$$
These upper bounds improve Mortici's inequality (\ref{ramamortici}).

\subsection{Volume of the unit ball}

The volume $\Omega_n$ of the unit ball in $\mathbb{R}^n$ is given in terms of the gamma function by the formula
$$
\Omega_n=\frac{\pi^{n/2}}{\Gamma(n/2+1)},\ \ n\in \mathbb{N}.
$$

Whereas the volume of the unit cube is $1$ in all dimensions, the numbers $\Omega_n$ strictly increase to the maximum $\Omega_5=8\pi^2/15$ and then strictly decrease  to $0$ as $n\to\infty$ (cf. \cite[p.264]{bh}).  G. D. Anderson, M. K. Vamanamurthy, and M. Vuorinen \cite{avv6} proved that $\Omega_n^{1/n}$ is strictly decreasing, and that the series $\sum_{n=2}^{\infty}\Omega_n^{1/\log n}$  is convergent.  In \cite{aq} G. D. Anderson and S.-L. Qiu proved that $\Omega_n^{1/(n\log n)}$ is strictly decreasing with limit $e^{-1/2}$ as $n\to \infty.$

 In 2008 H. Alzer published a collection of new inequalities for combinations of different dimensions and powers of $\Omega_n$ \cite[Section 3]{a4}.  We quote several of them below:
\begin{equation}\label{Alzer1}
a\frac{(2\pi e)^{n/2}}{n^{(n-1)/2}}\le (n+1)\Omega_n-n\Omega_{n+1}< b\frac{(2\pi e)^{n/2}}{n^{(n-1)/2}}, \quad n\ge1,
\end{equation}
where the best possible constants are $a=(4-9\pi/8)(2/(\pi e))^{1/2}/e=0.0829\dots$ and $b=\pi^{-1/2}=0.5641\ldots$;
\begin{equation}\label{Alzer3}
a\frac{(2\pi e)^n}{n^{n+2}}\le \Omega_n^2-\Omega_{n-1}\Omega_{n+1}<b\frac{(2\pi e)^n}{n^{n+2}},\quad n\ge2,
\end{equation}
with best possible constant factors $a=(4/e^2)(1-8/(3\pi))=0.0818\ldots$ and $b=1/(2\pi )=0.1591\ldots$;
\begin{equation}\label{Alzer4}
\frac{a}{\sqrt{n}}\le\frac{\Omega_n}{\Omega_{n-1}+\Omega_{n+1}}<\frac{b}{\sqrt{n}}, \quad n\ge2,
\end{equation}
with best possible constants $a=3\sqrt{2}\pi/(6+4\pi)=0.7178\ldots$ and $b=\sqrt{2\pi}=2.5066\ldots$; and
\begin{equation}\label{Alzer6}
\frac{a}{\sqrt{n}}\le (n+1)\frac{\Omega_{n+1}}{\Omega_n}-n\frac{\Omega_n}{\Omega_{n-1}}<\frac{b}{\sqrt{n}}, \quad n\ge2,
\end{equation}
with best possible constants $a=(4-\pi )\sqrt{2}=1.2139\ldots$ and $b=\sqrt{2\pi }/2=1.2533\ldots$.

Alzer's work in \cite{a4}  includes a number of new results about the gamma function and its derivatives.

In 2010 C. Mortici \cite{m2}, improving on some earlier work of H. Alzer \cite[Theorem 1]{a3}, obtained, for $n\ge 1$ on the left and for $n\ge 4$ on the right,
$$
\frac{a}{\sqrt[2n]{2\pi }}\le\frac{\Omega_n}{\Omega_{n+1}^{n/(n+1)}}<\frac{\sqrt{e}}{\sqrt[2n]{2\pi }},
$$
where $a=64\cdot 720^{11/12}\cdot 2^{1/22}/(10395\cdot \pi^{5/11})=1.5714\ldots.$  He sharpened work of H. Alzer \cite[Theorem 2]{a3} and S.-L. Qiu and M. Vuorinen \cite{qv} in the following result, valid for $n\ge 1$:
$$
\sqrt{\frac{2n+1}{4\pi}}<\frac{\Omega_{n-1}}{\Omega_n}<\sqrt{\frac{2n+1}{4\pi}+\frac{1}{16\pi n}}\ .
$$
C. Mortici also proved, in \cite[Theorem 4]{m2}, that, for $n\ge4,$
$$
\left(1+\frac{1}{n}\right)^{\frac{1}{2}-\frac{1}{4n}}<\frac{\Omega_n^2}{\Omega_{n-1}\Omega_{n+1}}<\left(1+\frac{1}{n}\right)^{\frac{1}{2}}.
$$
This result improves a similar one by H. Alzer \cite[Theorem 3, valid for $n\ge1$]{a3}, where the exponent on the left is the constant $2-\log_2\pi.$ Very recently, L. Yin \cite{yin} improved Mortici's result as follows: For $n\geq1$,
$$
\frac{(n+1)(n+1/6)}{(n+\beta)^2}<\frac{\Omega_n^2}{\Omega_{n-1}\Omega_{n+1}}<\frac{(n+1)(n+\beta/2)}{(n+1/3)^2},
$$
Where $\beta=(391/30)^{1/3}$.

\subsection{Digamma and polygamma functions}

The logarithmic derivative of the gamma function, $\psi(x)\equiv \frac{d}{dx}\log \Gamma(x) = \Gamma'(x)/\Gamma(x)$, is known as the \emph{digamma function}.  Its derivatives $\psi^{(n)}, n\ge 1,$ are known as the \emph{polygamma functions} $\psi_n.$  These functions have the following representations \cite[pp. 258--260]{as}  for $x>0$ and each natural number $n$:
\begin{equation*}
\psi(x)=-\gamma+\int_0^{\infty}\frac{e^{-t}-e^{-xt}}{1-e^{-t}}dt=-\gamma-\frac{1}{x}+\sum_{n=1}^{\infty}\frac{x}{n(x+n)}
\end{equation*}
and
\begin{equation*}
\psi_n(x)=(-1)^{n+1}\int_0^{\infty}\frac{t^ne^{-xt}}{1-e^{-t}}dt = (-1)^{n+1}n!\sum_{n=0}^{\infty}(x+k)^{-n-1} .
\end{equation*}

Several researchers have studied the properties of these functions.  In 2007, refining the left inequality in \cite[Theorem 4.8]{a0},  N. Batir \cite{bat3} obtained estimates for $\psi_n$ in terms of $\psi$ or $\psi_k$, with $k<n$.  In particular, he proved, for $x>0$ and $n\in\mathbb{N}$:
\begin{equation*}
(n-1)!\exp\left(-n\psi(x+1/2)\right)<\left|\psi_n(x)\right|<(n-1)!\exp\left(-n\psi(x)\right),
\end{equation*}
and, for $1\le k \le n-1$, $x>0$,
\begin{equation*}
(n-1)!\left(\frac{\psi_k(x+1/2)}{(-1)^{k-1}(k-1)!}\right)^{n/k}<|\psi_n(x)|<(n-1)!\left(\frac{\psi_k(x)}{(-1)^{k-1}(k-1)!}\right)^{n/k}.
\end{equation*}
He also proved, for example, the difference formula
\begin{equation*}
\alpha<\left((-1)^{n-1}\psi_n(x+1)\right)^{-1/n}-\left((-1)^{n-1}\psi_n(x)\right)^{-1/n}<\beta,
\end{equation*}
where $\alpha=(n!\zeta(n+1))^{-1/n}$ and $\beta=((n-1)!)^{-1/n}$ are best possible, and the sharp estimates
\begin{equation*}
-\gamma <\psi (x)+\log (e^{1/x}-1)<0,
\end{equation*}
where $\gamma$ is the Euler-Mascheroni constant.

In 2010 C. Mortici \cite{m3}  proved the following estimates, for $x>0$ and $n\ge 1$, refining work of B.-N. Guo, C.-P. Chen, and F. Qi \cite{gcq}:
$$
 -\frac{1}{720}\frac{(n+3)!}{x^{n+4}}<|\psi_n(x)|-\left[\frac{(n-1)!}{x^n}+\frac{1}{2}\frac{n!}{x^{n+1}}+\frac{1}{12}\frac{(n+1)!}{x^{n+2}}\right]<0.
$$

\subsection{Completely monotonic functions}\label{completemonot}

A function $f$ is said to be \emph{completely monotonic} on an interval $I$ if $(-1)^nf^{(n)}(x) \ge 0$ for all $x\in I$ and all nonnegative integers $n$.
If this inequality is strict, then $f$ is called \emph{strictly completely monotonic.}  Such functions occur in probability theory, numerical analysis, and other areas.  Some of the most important completely monotonic functions are the gamma function and the digamma and polygamma functions.  The Hausdorff-Bernstein-Widder theorem \cite[Theorem 12b, p. 161]{wi} states that  $f$ is completely monotonic on $[0,\infty)$ if and only if there is a non-negative measure $\mu$ on $[0,\infty)$ such that
$$
f(x)=\int_0^{\infty}e^{-xt}d\mu(t)
$$
for all $x>0$.
There is a well-written introduction to completely monotonic functions in \cite{ms}.

In 2008 N. Batir \cite{bat2} proved that the following function $F_a(x)$  related to the gamma function is completely monotonic on $(0,\infty)$ if and only if $a\ge 1/4$ and that $-F_a(x)$ is completely monotonic if and only if $a\le0$:
\begin{equation*}
F_a(x)\equiv \log\Gamma(x)-x\log x + x -\frac{1}{2}\log(2\pi)+\frac{1}{2}\psi(x)+\frac{1}{6(x-a)}.
\end{equation*}
As a corollary he was able to prove, for $x>0$, the inequality
\begin{equation*}
\exp\left(-\frac{1}{2}\psi(x)-\frac{1}{6(x-\alpha)}\right)<\frac{\Gamma(x)}{x^xe^{-x}\sqrt{2\pi}}<\exp\left(-\frac{1}{2}\psi(x)-\frac{1}{6(x-\beta)}\right),
\end{equation*}
with best constants $\alpha=1/4$ and $\beta=0,$ improving his earlier work with H. Alzer \cite{alb}.

In 2010 C. Mortici \cite{m3} showed that for every $n\ge 1$, the functions $f,g:(0,\infty)\to \mathbb{R}$ given by
\begin{equation*}
f(x)\equiv |\psi_n(x)|-\frac{(n-1)!}{x^n}-\frac{1}{2}\frac{n!}{x^{n+1}}-\frac{1}{12}\frac{(n+1)!}{x^{n+2}}+\frac{1}{720}\frac{(n+3)!}{x^{n+4}}
\end{equation*}
and
\begin{equation*}
g(x)\equiv \frac{(n-1)!}{x^n}+\frac{1}{2}\frac{n!}{x^{n+1}}+\frac{1}{12}\frac{(n+1)!}{x^{n+2}}-|\psi_n(x)|
\end{equation*}
are completely monotonic on $(0,\infty)$.  As a corollary, since  $f(x)$ and $g(x)$ are positive, he obtained estimates for $|\psi_n(x)|$ as finite series in negative powers of $x$.

G. D. Anderson and S.-L. Qiu \cite{aq}, as well as some other authors (see \cite{aa}), have studied the monotonicity properties of the function   $f(x)\equiv (\log\Gamma(x+1))/x$.  In 2011 J. A. Adell and H. Alzer \cite{aa} proved that  $f'$ is completely monotonic on $(-1,\infty).$

In the course of pursuing research inspired by \cite{aq} and \cite{avv6}  (see \cite{bp}), in 2012 H. Alzer \cite{a4} discussed properties of the function
\begin{equation*}
f(x)\equiv \left(1-\frac{\log x}{\log (1+x)}\right) x\log x,
\end{equation*}
which F. Qi and B.-N. Guo \cite{qg} later conjectured to be completely monotonic on $(0,\infty)$.  In \cite{bp} C. Berg and H. L. Pedersen  proved this conjecture.

In 2001 C. Berg and H. L. Pedersen \cite{bp1} proved that the derivative of the function
$$f(x)\equiv\frac{\log\Gamma(x+1)}{x\log x},\quad x>0,$$
is completely monotonic (see also \cite{bp2}). This result extends work of \cite{aq} and \cite{el}.
Very recently, C. Berg and H. L. Pedersen \cite{bp3} showed that the function
$$F_a(x)\equiv\frac{\log\Gamma(x+1)}{x\log(ax)}$$
is a Pick function when $a\geq1$, that is, it extends to a holomorphic function mapping the upper half plane into itself. The authors also considered the function
$$f(x)\equiv\left(\frac{\pi^{x/2}}{\Gamma(1+x/2)}\right)^{1/(x\log x)}$$
and proved that $\log f(x+1)$ is a Stieltjes function and hence that $f(x+1)$ is completely monotonic on $(0,\infty)$.

\section{The hypergeometric function and Elliptic integrals}

The classical \emph{hypergeometric function} is defined by
$$
F(a,b;c;x)\equiv {}_2F_1(a,b;c;x)=\sum_{n=0}^{\infty}\frac{(a,n)(b,n)}{(c,n)}\frac{x^n}{n!}, \ \ |x|<1,
$$
where $(a,n)\equiv a(a+1)(a+2)\cdots (a+n-1)$ for $n\in \mathbb{N}$, and $(a,0)=1$ for $a\ne 0$.  This function is so general that for proper choice of the parameters $a$, $b$, $c$, one obtains logarithms, trigonometric functions, inverse trigonometric functions, elliptic integrals, or polynomials of Chebyshev, Legendre, Gegenbauer, or Jacobi, and so on (see \cite[ch. 15]{as}).

\subsection{Hypergeometric functions}

The Bernoulli inequality \cite[p. 34]{mi} may be written as
\begin{equation}\label{bernoulli}
\log (1+ct)\le c\log (1+t),
\end{equation}
where  $c>1$, $t>0$. In \cite{kmv} some Bernoulli-type inequalities have been obtained.

It is well known that in the zero-balanced case $c=a+b$ the hypergeometric function $F(a,b;c;x)$ has a logarithmic singularity at $x=1$ (cf. \cite[Theorem 1.19(6)]{avv2}).  Moreover, as a special case \cite[15.1.3]{as},
\begin{equation}\label{log}
xF(1,1,2;x)=\log\frac{1}{1-x}.
\end{equation}

Because of this connection, M. Vuorinen and his collaborators \cite{kmsv} have generalized  versions of \eqref{bernoulli} to a wide class of hypergeometric functions.  In the course of their investigation they have studied monotonicity and convexity/concavity properties of such functions.  For example, for positive $a, b$ let $g(x)\equiv xF(a,b;a+b;x)$, $x\in (0,1).$ These authors have proved that $G(x)\equiv\log g(e^x/(1+e^x))$ is concave on $(-\infty,\infty)$ if and only if $1/a+1/b\ge 1.$  And they have shown that, for fixed $a, b\in (0,1]$ and for $ x\in (0,1)$, $p>0$, the function
$$
\left(\frac{x^p}{1+x^p}F\left(a,b;a+b;\frac{x^p}{1+x^p}\right)\right)^{1/p}
$$
is increasing in $p$. In particular,
$$
\frac{\sqrt r}{1+\sqrt r}F\left(a,b;a+b;\frac{\sqrt r}{1+\sqrt r} \right)\leq\left(\frac{r}{1+r}F\left(a,b;a+b;\frac{r}{1+r}\right)\right)^{1/2}.
$$

Motivated by the asymptotic behavior of $F(x)=F(a,b;c;x)$ as $x\to 1^{-}$, S. Simi\'c and M. Vuorinen have carried the above work further in \cite{sv2}, finding best possible bounds, when $a,b,c>0$ and $0<x, y<1$,  for the quotient and difference
$$
\frac{F(a,b;c;x)+F(a,b;c;y)}{F(a,b;c;x+y-xy)},\ \ \ F(x)+F(y)-F(x+y-xy) .
$$

In 2009 D. Karp and S. M. Sitnik \cite{ks} obtained some inequalities and monotonicity of ratios for generalized hypergeometric function. The proofs hinge on a generalized Stieljes representation of the generalized hypergeometric function.

\subsection{Complete elliptic integrals}

For $0<r<1$, the \emph{complete elliptic integrals of the first and second kind} are defined as
\begin{equation}\label{ellipK}
\ellipK(r) \equiv\int_0^{\pi/2}\frac{dt}{\sqrt{1-r^2\sin^2t}}=\int_0^1\frac{dt}{\sqrt{(1-t^2)(1-r^2t^2)}}
\end{equation}
and
\begin{equation}\label{ellipE}
\ellipE(r)
\equiv\int_0^{\pi/2}{\sqrt{1-r^2\sin^2t}}\, dt=\int_0^1\sqrt{\frac{1-r^2t^2}{1-t^2}}\, dt,
\end{equation}
respectively.
Letting $r'\equiv\sqrt{1-r^2}$, we often denote
$$\ellipK'(r) =\ellipK(r'),\quad \ellipE'(r) =\ellipE(r').$$
These elliptic integrals have the hypergeometric series representations
\begin{equation}\label{ellip2F1}
\mathcal{K}(r)=\frac{\pi}{2}F\left(\tfrac{1}{2},\tfrac{1}{2};1;r^2\right),\ \ \ \mathcal{E}=\frac{\pi}{2}F\left(\tfrac{1}{2},-\tfrac{1}{2};1;r^2\right)\ .
\end{equation}

\subsection{The Landen identities}

The functions $\ellipK$ and $\ellipE$ satisfy the following identities due to Landen \cite[163.01, 164.02]{bf}:
$$
\ellipK\left(\dfrac{2\sqrt{r}}{1+r}\right)=(1+r)\ellipK(r) ,
\quad
\ellipK\left(\dfrac{1-r}{1+r}\right)=\dfrac12(1+r)\ellipK'(r) ,
$$
$$
\ellipE\left(\dfrac{2\sqrt{r}}{1+r}\right)=\dfrac{2\ellipE(r) -r'^2\ellipK(r) }{1+r},
\quad
\ellipE\left(\dfrac{1-r}{1+r}\right)=\dfrac{\ellipE'(r) +r\ellipK'(r) }{1+r}.
$$

Using Landen's transformation formulas, we have the following identities \cite[Lemma 2.8]{vz}:
For $r\in(0,1)$, let $t=(1-r)/(1+r)$. Then
$$
\ellipK(t^2)=\dfrac{(1+r)^2}4\ellipK'(r^2),
\quad
\ellipK'(t^2)={(1+r)^2}\ellipK(r^2),
$$
$$
\ellipE(t^2)=\dfrac{\ellipE'(r^2)+(r+r^2+r^3)\ellipK'(r^2)}{(1+r)^2},
$$
$$
\ellipE'(t^2)=\dfrac{4\ellipE(r^2)-(3-2r^2-r^4)\ellipK(r^2)}{(1+r)^2}.
$$

Generalizing a Landen identity, S. Simi\'c and M. Vuorinen \cite{sv1} have determined the precise regions in the $ab$-plane for which a Landen inequality holds for zero-balanced hypergeometric functions.  They proved that for all $a, b > 0$ with $ab\le 1/4$ the inequality
$$
F\left(a,b;a+b;\frac{4r}{(1+r)^2}\right)\le (1+r)F\left(a,b;a+b;r^2\right)
$$
holds for $r\in (0,1)$, while for $a, b >0$ with $1/a+1/b\le 4,$ the following reversed inequality is true for each $r\in(0,1):$
$$
F\left(a,b;a+b;\frac{4r}{(1+r)^2}\right)\ge (1+r)F\left(a,b;a+b;r^2\right)\ .
$$
In the rest of the $ab$-plane neither of these inequalities holds for all $r\in (0,1).$  These authors have also obtained sharp bounds for the quotient
$$
\frac{(1+r)F(a,b;a+b;r^2)}{F(a,b;a+b;4r/(1+r)^2)}
$$
in certain regions of the $ab$-plane.

Some earlier results on Landen inequalities for hypergeometric functions can be found in \cite{qv2}. Recently, \'A. Baricz obtained Landen-type inequalities for generalized Bessel functions \cite{bar3, bar4}.

Inspired by an idea of Simi\'c and  Vuorinen \cite{sv1}, M.-K. Wang, Y.-M. Chu, and Y.-P. Jiang \cite{wcj} obtained some inequalities for zero-balanced hypergeometric functions which generalize Ramanujan's cubic transformation formulas.

\subsection{Legendre's relation and generalizations}
It is well known that the complete elliptic integrals satisfy the \emph{Legendre relation} \cite[110.10]{bf}:
$$\ellipE\ellipK'+\ellipE'\ellipK-\ellipK\ellipK'=\frac{\pi}{2}.$$  		
This relation has been generalized in various ways.  E. B. Elliott \cite{e} proved the identity
$$
F_1F_2+F_3F_4-F_3F_2=\frac{\Gamma(1+\lambda+\mu)\Gamma(1+\mu+\nu)}{\Gamma(\lambda+\mu+\nu+\tfrac{3}{2})\Gamma(\tfrac{1}{2}+\mu)},
$$
where
$$F_1=F(\tfrac{1}{2}+\lambda,-\tfrac{1}{2}-\nu;1+\lambda+\mu;x),\ \ F_2=F(\tfrac{1}{2}-\lambda,\tfrac{1}{2}+\nu;1+\mu+\nu;1-x),$$
$$F_3=F(\tfrac{1}{2}+\lambda,\tfrac{1}{2}-\nu;1+\lambda+\mu;x),\ \ F_4=F(-\tfrac{1}{2}-\lambda,\tfrac{1}{2}+\nu;1+\mu+\nu;1-x).$$
Elliott proved this formula by a clever change of variables in multiple integrals.  Another proof, based on properties of the hypergeometric differential equation, was suggested without details in \cite[p.~138]{aar}, and the missing details were provided in \cite{avv3}.  It is easy to see that Elliott's formula reduces to the Legendre relation when  $\lambda=\mu=\nu=0$ and $x=r^2$.

Another generalization of the Legendre relation was given in \cite{aqvv}.  With the notation
$$u=u(r)=F(a-1,b;c;r),\ v=v(r)=F(a,b;c;r),$$
$$u_1=u(1-r),\ v_1=v(1-r),$$	
the authors considered the function
$$\mathcal{L}(a,b,c,r)=uv_1+u_1v-vv_1,$$
proving, in particular, that
$$\mathcal{L}(a,1-a,c,r)=\frac{\Gamma^2(c)}{\Gamma(c+a-1)\Gamma(c-a+1)}.$$
This reduces to Elliott's formula in case $\lambda=\nu=1/2-a$ and  $\mu=c+a-3/2$. In \cite{aqvv} it was conjectured that for $a, b\in(0,1)$, $a+b\leq1(\geq1)$, $\mathcal{L}(a,b,c,r)$ is concave (convex) as a function of  $r$  on  $(0,1)$.  In \cite{kv} Karatsuba and Vuorinen determined, in particular, the exact regions of $abc-$space in which the function  $\mathcal{L}(a,b,c,r)$  is concave, convex, constant, positive, negative, zero, and where it attains its unique extremum.

In \cite{bnpv} Balasubramanian, Naik, Ponnusamy, and Vuorinen obtained a differentiation formula for an expression involving hypergeometric series that implies Elliott's identity.  This paper contains a number of other significant results, including a proof that Elliott's identity is equivalent to a formula of Ramanujan \cite[p.~87, Entry 30]{bt} on the differentiation of quotients of hypergeometric functions.

\subsection{Some approximations for $\ellipK(r) $ by $\arth r$}
G. D. Anderson, M. K. Vamanamurthy, and M. Vuorinen \cite{avv7} approximated $\ellipK(r) $ by the inverse hyperbolic tangent function $\arth$, obtaining the inequalities
\begin{equation}\label{approxK:avv}
\frac{\pi}{2}\left(\frac{\arth r}{r}\right)^{1/2}<\ellipK(r) <\frac{\pi}{2}\frac{\arth r}{r},
\end{equation}
for $0<r<1$. H. Alzer and S.-L. Qiu \cite{alq} refined (\ref{approxK:avv}) as
\begin{equation}\label{approxK:alq}
\frac{\pi}{2}\left(\frac{\arth r}{r}\right)^{3/4}<\ellipK(r) <\frac{\pi}{2}\frac{\arth r}{r},
\end{equation}
with the best exponents $3/4$ and $1$ for $(\arth r)/r$ on the left and right, respectively.
Seeking to improve the exponents in \eqref{approxK:alq}, they conjectured that
the double inequality
\begin{equation}\label{approxK:alqcoj}
\frac{\pi}{2}\left(\frac{\arth r}{r}\right)^{3/4+\alpha r}<\ellipK(r) <\frac{\pi}{2}\left(\frac{\arth r}{r}\right)^{3/4+\beta r}
\end{equation}
holds for all $0<r<1$, with best constants $\alpha=0$ and $\beta=1/4$. Very recently Y.-M. Chu et al \cite{cwq2} gave a proof for this conjecture.

S. Andr\'as and \'A. Baricz \cite{ab} presented some improved lower and upper bounds for $\ellipK(r) $ involving the Gaussian hypergeometric series.

\subsection{Approximations for $\mathcal{E}(r)$}
In \cite{gq} B.-N. Guo and F. Qi have obtained new approximations for $\mathcal{E}(r)$ as well as for $\mathcal{K}(r)$.  For example, they showed that, for $0<r<1$,
$$
\frac{\pi}{2}-\frac{1}{2}\log\frac{(1+r)^{1-r}}{(1-r)^{1+r}}<\mathcal{E}(r)<\frac{\pi -1}{2}+\frac{1-r^2}{4r}\log \frac{1+r}{1-r}.
$$
In recent work \cite{cwqj, wc, wcqj} Y.-M. Chu et al have obtained estimates for $\mathcal{E}(r)$ in terms of rational functions of the arithmetic, geometric, and root-square means, implying new inequalities for the perimeter of an ellipse.

\subsection{Generalized complete elliptic integrals}

For $0<a<\min\{c,1\}$ and $0<b<c\le a+b$, define the \emph{generalized complete elliptic integrals of the first and second kind} on $[0,1]$ by \cite{hvv}
\begin{equation}\label{eqn:GenEllipK}
\ellipK_{a,b,c}=\ellipK_{a,b,c}(r)\equiv\frac{B(a,b)}{2}F(a,b;c;r^2),
\end{equation}
\begin{equation}\label{eqn:GenEllipE}
\ellipE_{a,b,c}=\ellipE_{a,b,c}(r)\equiv\frac{B(a,b)}{2}F(a-1,b;c;r^2),
\end{equation}
\begin{equation}\label{eqn:GenEllipKEprime}
\ellipK'_{a,b,c}=\ellipK_{a,b,c}(r')\quad\mbox{and}\quad\,\ellipE'_{a,b,c}=\ellipE_{a,b,c}(r'),
\end{equation}
for $r\in(0,1)$, $r'=\sqrt{1-r^2}$. The end values are defined by limits as $r$ tends to $0+$ and $1-$, respectively. Thus,
$$
\ellipK_{a,b,c}(0)=\ellipE_{a,b,c}(0)=\frac{B(a,b)}{2}
$$
and
$$
\ellipE_{a,b,c}(1)=\frac12\frac{B(a,b)B(c,c+1-a-b)}{B(c+1-a,c-b)},\quad \ellipK_{a,b,c}(1)=\infty.
$$
Note that the restrictions on the parameters $a$, $b$, and $c$ ensure that the function $\ellipK_{a,b,c}$ is increasing and unbounded, whereas $\ellipE_{a,b,c}$ is decreasing and bounded, as in the classical case $a=b=1/2$, $c=1$.

V. Heikkala, H. Lind\'en, M. K. Vamanamurthy, and M. Vuorinen \cite{hvv,hlvv} derived several differentiation formulas, and obtained sharp monotonicity and convexity properties for certain combinations of the generalized elliptic integrals. They also constructed a conformal mapping $\rm{sn}_{a,b,c}$ from a quadrilateral with internal angles $b\pi$, $(c-b)\pi$, $(1-a)\pi$, and $(1-c+a)\pi$ onto the upper half plane. These results generalize the work of \cite{aqvv}. For some particular parameter triples $(a,b,c)$ there are very recent results by many authors \cite{bar4, wzc1, zwc2, zqw}.

With suitable restrictions on the parameters $a,b,c$, E. Neuman \cite{n2} has obtained bounds for $\mathcal{K}_{a,b,c}$ and $\mathcal{E}_{a,b,c}$ and for certain combinations and products of them.  He has also proved that these generalized elliptic integrals are logarithmically convex as functions of the first parameter.

In 2007 \'A. Baricz \cite{bar7, bar8, bar10} established some Tur\'an-type inequalities for Gaussian hypergeometric functions and generalized complete elliptic integrals. He also studied the generalized convexity of the zero-balanced hypergeometric functions and generalized complete elliptic integrals \cite{bar2} (see also \cite{bar5, bar9, bar4}). Very recently, S.I. Kalmykov and D.B. Karp \cite{kk1, kk2} have studied log-convexity and log-concavity for series involving gamma functions and derived many known and new inequalities for the modified Bessel, Kummer and generalized hypergeometric functions and ratios of the Gauss hypergeometric functions. In particular, they improved and generalized Baricz's Tur\'an-type inequalities.

\subsection{The generalized modular function and generalized linear distortion function}
Let $a,b,c>0$ with $a+b\geq c$.
A \emph{generalized modular equation of order} (or \emph{degree}) $p>0$ is
\begin{equation}\label{eqn:modular}
\frac{F(a,b;c;1-s^2)}{F(a,b;c;s^2)}=p\frac{F(a,b;c;1-r^2)}{F(a,b;c;r^2)},\quad\,0<r<1.
\end{equation}
The \emph{generalized modulus} is the decreasing homeomorphism $\mu_{a,b,c}:(0,1)\to(0,\infty)$, defined by
\begin{equation}\label{eqn:mu}
\mu_{a,b,c}(r)\equiv\frac{B(a,b)}{2}\frac{F(a,b;c;1-r^2)}{F(a,b;c;r^2)}.
\end{equation}
The generalized modular equation (\ref{eqn:modular}) can be written as
$$
\mu_{a,b,c}(s)=p\mu_{a,b,c}(r).
$$
With $p=1/K$, $K>0$, the solution of (\ref{eqn:modular}) is then given by
\begin{equation*}
s=\varphi_K^{a,b,c}(r)\equiv\mu_{a,b,c}^{-1}(\mu_{a,b,c}(r)/K).
\end{equation*}
Here $\varphi_{K}^{a,b,c}$ is called the $(a,b,c)$-\emph{modular function} with \emph{degree} $p=1/K$ \cite{aqvv, hlvv, hvv}.
Clearly the following identities hold:
$$
\mu_{a,b,c}(r)\mu_{a,b,c}(r')=\left(\frac{B(a,b)}{2}\right)^2,
$$
$$
\varphi_K^{a,b,c}(r)^2+\varphi_K^{a,b,c}(r')^2=1.
$$
In \cite{hlvv}, the authors generalized the functional inequalities for the modular functions and Gr\"otzsch function $\mu$ proved in \cite{aqvv} to hold also for the generalized modular functions and generalized modulus in the case $b=c-a$.
For instance, for $0<a<c\leq1$ and  $K>1$, the inequalities
\begin{equation}\label{mu-chain-inequal}
\mu_{a,c-a,c}(1-\sqrt{(1-u)(1-t)})\leq\frac{\mu_{a,c-a,c}(u)+\mu_{a,c-a,c}(t)}{2}\leq\mu_{a,c-a,c}(\sqrt{ut})
\end{equation}
hold for all $u,t\in(0,1)$, with equality if and only if $u=t$, and
\begin{equation}
r^{1/K}<\varphi_K^{a,c-a,c}(r)<e^{(1-1/K)R(a,c-a)/2}r^{1/K},
\end{equation}
\begin{equation}
r^{K}>\varphi_K^{a,c-a,c}(r)>e^{(1-K)R(a,c-a)/2}r^{K}.
\end{equation}
For the special case of $a=1/2$ and $c=1$ the readers are referred to \cite{avv2}.
G.-D. Wang et al \cite{wzc1} presented several sharp inequalities for the generalized modular functions with specific choice of parameters $c=1$ and $b=1-a$.

A linearization for the generalized modular function is also presented in \cite{hlvv} as follows:
Let $p:(0,1)\to(-\infty,\infty)$ and $q:(-\infty,\infty)\to(0,1)$ be given by $p(x)=2\log(x/x')$ and $q(x)=p^{-1}(x)=\sqrt{e^x/(e^x+1)}$, respectively, and for $a\in(0,1)$, $c\in(a,1]$, $K\in(1,\infty)$, let $g,h:(-\infty,\infty)\to(-\infty,\infty)$ be defined by $g(x)=p(\varphi_K^{a,c-a,c}(q(x)))$ and $h(x)=p(\varphi_{1/K}^{a,c-a,c}(q(x)))$. Then
$$
g(x)\geq\left\{
\begin{array}{ll}
Kx,&{\rm if}\quad x\geq0,\\
x/K,&{\rm if}\quad x<0,
\end{array}\right.\quad{\rm and}\quad
h(x)\leq\left\{
\begin{array}{ll}
x/K,&{\rm if}\quad x\geq0,\\
Kx,&{\rm if}\quad x<0.
\end{array}\right.
$$
In the same paper the authors also studied how these generalized functions depend on the parameter $c$.  Corresponding results for the case $c=1$ can be found in the articles \cite{aqvv, qv3, zwc3}.

Recently B. A. Bhayo and M. Vuorinen \cite{bv} have studied the H\"older continuity and submultiplicative properties of $\varphi_K^{a,b,c}(r)$ in the case where $c=1$ and $b=1-a$, and have obtained several sharp inequalities for $\varphi_K^{a,1-a,1}(r).$

For $x, K\in (0,\infty ),$ define
$$
\eta_K^a(x)\equiv\left(\frac{s}{s'}\right)^2,\ s=\varphi_K^{a,1-a,1}(r),\ r=\sqrt{\frac{x}{1+x}},
$$
and the \emph{generalized linear distortion function}
$$
\lambda (a,K)\equiv \left(\frac{\varphi_K^{a,1-a,1}(1/\sqrt{2})}{\varphi_{1/K}^{a,1-a,1}(1/\sqrt{2})}\right)^2,\ \ \lambda (a,1)=1.
$$
For $a=1/2$, these two functions reduce to the well-known special case denoted by  $\eta_K(x)$ and $\lambda(K)$, respectively,  which play a crucial role in quasiconformal theory.  Several inequalities for these functions have been obtained as an application of the monotonicity and convexity of certain combinations of these functions and some elementary functions; see \cite{bv, cwzq, mqzc, mqzc2, wzqc, zwcq}. For instance, the following chain of inequalities is proved in \cite{cwzq}: for $a\in(0,1/2]$, $K\in(1,\infty)$, and $x,y\in(0,\infty)$,
$$
\max\left\{\frac{2\eta_K^a(x)\eta_K^a(y)}{\eta_K^a(x)+\eta_K^a(y)},\eta_K^a\left(\frac{2xy}{x+y}\right)\right\}
\leq\eta_K^a(\sqrt{xy})\hspace{1cm}
$$
$$
\hspace{1cm}\leq\sqrt{\eta_K^a(x)\eta_K^a(y)}\leq\min\left\{\frac{\eta_K^a(x)+\eta_K^a(y)}{2},\eta_K^a\left(\frac{x+y}{2}\right)\right\},
$$
with equality if and only if $x=y$.

\section{Inequalities for power series}

The following theorem is an interesting tool in simplified proofs for monotonicity of the quotient of two power series.

\begin{theorem}\cite[Theorem 4.3]{hvv}\label{rule4series}
Let $\sum_{n=0}^\infty a_n x^n$ and $\sum_{n=0}^\infty b_n x^n$ be two real power series converging on the interval $(-R,R)$. If the sequence $\{a_n/b_n\}$ is increasing (decreasing) and $b_n>0$ for all $n$, then the function
$$
f(x)=\frac{\sum_{n=0}^\infty a_n x^n}{\sum_{n=0}^\infty b_n x^n}
$$
is also increasing (decreasing) on $(0,R)$. In fact, the function
$$f'(x)(\sum_{n=0}^\infty b_n x^n)^2$$
has positive Maclaurin coefficients.
 \end{theorem}

A more general version of this theorem appears in \cite{bk} and \cite[Lemma 2.1]{pv}. This kind of rule also holds for the quotient of two polynomials instead of two power series.

\begin{theorem}\cite[Theorem 4.4]{hvv}
Let $f_n(x)=\sum_{k=0}^n a_k x^k$ and $g_n(x)=\sum_{k=0}^n b_k x^k$ be two real polynomials, with $b_k>0$ for all $k$. If the sequence $\{a_k/b_k\}$ is increasing (decreasing), then so is the function $f_n(x)/g_n(x)$ for all $x>0$. In fact, $g_n f_n'-f_n g_n'$ has positive (negative) coefficients.
\end{theorem}

In 1997 S. Ponnusamy and M. Vuorinen \cite{pv} refined Ramanujan's work on asymptotic behavior of the hypergeometric function and also obtained many inequalities for the hypergeometric function by making use of Theorem \ref{rule4series}.
Many well-known results of monotonicity and inequalities for complete elliptic integrals have been extended to the generalized elliptic integrals in \cite{hlvv, hvv}.

Motivated by an open problem of G. D. Anderson et al \cite{avv9}, in 2006 \'A. Baricz \cite{bar5} considered ratios of general power series and obtained the following theorem. Note the similarity of the last inequality in Theorem \ref{baricz} with the left-hand side of the inequality (\ref{mu-chain-inequal}).

\begin{theorem}\label{baricz}
Suppose that the power series $f(x)=\sum_{n=0}^\infty a_n x^n$ with $a_n>0$ for all $n\geq0$ is convergent for all $x\in(0,1)$, and also that the sequence $\{(n+1)a_{n+1}/a_n-n\}_{n\geq0}$ is strictly decreasing. Let the function $m_f:(0,1)\to(0,\infty)$ be defined as $m_f(r)=f(1-r^2)/f(r^2)$. Then
$$
\sqrt[\uproot{8}k]{\prod_{i=1}^k m_f(r_i)}\leq m_f\left(\sqrt[\uproot{8}k]{\prod_{i=1}^k r_i}\right),
$$
for all $r_1,r_2,\cdots,r_k\in(0,1)$, where equality holds if and only if $r_1=r_2=\cdots=r_k$. In particular, for $k=2$ the inequalities
$$
\sqrt{m_f(r_1)m_f(r_2)}\leq m_f(\sqrt{r_1 r_2}),
$$
$$
\frac{1}{m_f(r_1)}+\frac{1}{m_f(r_2)}\geq\frac{2}{m_f(\sqrt{r_1 r_2})},
$$
$$
m_f(r_1)+m_f(r_2)\geq2m_f\left(\sqrt{1-\sqrt{(1-r_1^2)(1-r_2^2)}}\right)
$$
hold for all $r_1,r_2\in(0,1)$, and in all these inequalities equality holds if and only if $r_1=r_2$.
\end{theorem}

The following Landen-type inequality for power series is also due to Baricz \cite{bar3}.

\begin{theorem}
Suppose that the power series $f(x)=\sum_{n=0}^\infty a_n x^n$ with $a_n>0$ for all $n\geq0$ is convergent for all $x\in(0,1)$, and that for a given $\delta>1$ the sequence $\{n!a_n/(\log\delta)^n\}_{n\geq0}$ is decreasing. If $\lambda_f(x)=f(x^2)$, then
$$
\lambda_f\left(\frac{2\sqrt{r}}{1+r}\right)<\rho \lambda_f(r)
$$
holds for all $r\in(0,1)$ and $\rho\geq\delta^{4\sqrt{2}-5}$.
\end{theorem}

G. D. Anderson, M. K. Vamanamurthy, and M. Vuorinen \cite{avv8} studied generalized convexity and gave sufficient conditions for generalized convexity of functions defined by Maclaurin series. These results yield a class of new inequalities for power series which improve some earlier results obtained by \'A. Baricz. More inequalities for power series can be found in \cite{bar4, cwzq}.

In 1928 T. Kaluza gave a criterion for the signs of the power series of a function that is the reciprocal of another power series.

\begin{theorem}\cite{kal}
Let $f(x)=\sum_{n\geq0}a_n x^n$ be a convergent Maclaurin series with radius of convergence $r>0$. If $a_n>0$ for all $n\geq0$ and the sequence $\{a_n\}_{n\geq0}$ is log-convex, that is, for all $n\geq0$
\begin{equation}\label{log-condition}
a_n^2\leq a_{n-1}a_{n+1},
\end{equation}
then the coefficients $b_n$ of the reciprocal power series $1/f(x)=\sum_{n\geq0}b_n x^n$ have the following properties: $b_0=1/a_0>0$ and $b_n\leq0$ for all $n\geq1$.
\end{theorem}
In 2011 \'A. Baricz, J. Vesti, and M. Vuorinen \cite{bvv} showed that the condition (\ref{log-condition}) cannot be replaced by the condition
$$
a_n\leq\left(\frac{a_{n-1}^t+a_{n+1}^t}{2}\right)^{1/t},
$$
for any $t>0$. However, it is not known whether the condition (\ref{log-condition}) is necessary.

In 2009 S. Koumandos and H. L. Pedersen \cite{kp} proved the following interesting result, which deals with the monotonicity properties of the quotient of two series of functions.

\begin{theorem}\cite[Lemma 2.2]{kp}
Suppose that $a_k>0$, $b_k>0$ and that $\{u_k(x)\}$ is a sequence of positive $C^1$-functions such that the series
$$
\sum_{k=0}^\infty a_k u_k^{(l)}(x)\quad {\rm and }\quad \sum_{k=0}^\infty b_k u_k^{(l)}(x),\quad l=0,\, 1,
$$
converge absolutely and uniformly over compact subsets of $[0,\infty)$. Define
$$
f(x)\equiv\frac{\sum_{k=0}^\infty a_k u_k(x)}{\sum_{k=0}^\infty b_k u_k(x)}.
$$
(1) If the logarithmic derivatives $u_k'(x)/u_k(x)$ form an increasing sequence of functions and if $a_k/b_k$ decreases (resp. increases) then $f(x)$ decreases (resp. increases) for $x\geq0$.\\
(2) If the logarithmic derivatives $u_k'(x)/u_k(x)$ form a decreasing sequence of functions and if $a_k/b_k$ decreases (resp. increases) then $f(x)$ increases (resp. decreases) for $x\geq0$.
\end{theorem}

For inequalities of power series as complex functions, see \cite{id, idcd, idd} and the references therein.

\section{Means}

A \emph{homogeneous bivariate mean} is defined as a continuous function $\mathcal{M}:\R^+\times\R^+\to\R$
satisfying $\min\{x,y\}\leq\mathcal{M}(x,y)\leq\max\{x,y\}$ and $\mathcal{M}(\lambda x,\lambda y)=\lambda\mathcal{M}(x,y)$ for all $x,y,\lambda>0$.
Important examples are the \emph{arithmetic mean} $A(a,b)$, the \emph{ geometric mean} $G(a,b)$, the \emph{logarithmic mean} $L(a,b)$, the \emph{identric mean} $I(a,b)$, the \emph{root square mean} $Q(a,b)$, and the \emph{power mean} $M_r(a,b)$ \emph{of order} $r$, defined, respectively, by
$$
A(a,b)=\frac{a+b}{2},\ \ G(a,b)=\sqrt{ab},
$$
$$
L(a,b)=\frac{a-b}{\log a-\log b},\ \ I(a,b)=\frac{1}{e}\left(\frac{a^a}{b^b}\right)^{1/(a-b)},
$$
$$
Q(a,b)=\sqrt{\frac{a^2+b^2}{2}},\ \ M_r(a,b)=\sqrt[\leftroot{-3}\uproot{8}r]{\frac{a^r+b^r}{2}}.
$$

\subsection{Power means}
The \emph{weighted power means} are defined by
\begin{equation*}
M_\lambda(\omega;a,b)\equiv\left(\omega a^\lambda+(1-\omega)b^\lambda\right)^{1/\lambda}\qquad (\lambda\neq0),
\end{equation*}
$M_0(\omega;a,b)\equiv a^\omega b^{1-\omega}$, with \emph{weights} $\omega,1-\omega>0$. The \emph{power means} are the equally-weighted means $M_\lambda(a,b)=M_\lambda(1/2;a,b)$. As a special case, we have $M_0(1/2;a,b)=G(a,b)$.

In \cite{k2} O. Kouba studied the ratio of differences of power means
$$
\rho (s,t,p;a,b) \equiv \frac{M_s^p(a,b)-G^p(a,b)}{M_t^p(a,b)-G^p(a,b)},
$$
finding sharp bounds for $\rho (s,t,p;a,b)$ in various regions of $stp$-space with $a,b$ positive and $a\ne b$.  This work extends results of H. Alzer and S.-L. Qiu \cite{alq1}, T. Trif \cite{t}, O. Kouba \cite{k1}, S.-H. Wu \cite{wu}, and S.-H. Wu and L. Debnath \cite{wd3}.  O. Kouba also extended the range of validity of the following inequality, due to S.-H. Wu and L. Debnath \cite{wd3}:
$$
\frac{2^{-p/r}-2^{-p/s}}{2^{-p/t}-2^{-p/s}}<\frac{M_r^p(a,b)-M_s^p(a,b)}{M_t^p(a,b)-M_s^p(a,b)}< \frac{r-s}{t-s}
$$
to the set of real numbers $r,t,s,p$ satisfying the conditions $0<s<t<r$ and $0<p\le (4t+2s)/3.$

\subsection{Toader means}

If $p:\R^+\to\R^+$ is a strictly monotonic function, then define
$$
f(a,b;p,n)\equiv\left\{
\begin{array}{ll}
\dfrac{1}{2\pi}\int_0^{2\pi}p((a^n\cos^2\theta+b^n\sin^2\theta)^{1/n})d\theta & {\rm if}\quad n\neq0,\vspace{3pt}\\
\dfrac{1}{2\pi}\int_0^{2\pi}p(a^{\cos^2\theta}b^{\sin^2\theta})d\theta & {\rm if}\quad n=0,
\end{array}
\right.
$$
where $a,b$ are positive real numbers. The \emph{Toader mean} \cite{to} of $a$ and $b$ is defined as $T(a,b;p,n)\equiv  p^{-1}(f(a,b;p,n))$. It is easy to see that the Toader mean is symmetric. For special choices of $p$, let
$T_{q,n}(a,b)=T(a,b;p,n)$ if $p(x)=x^q$ with $q\neq0$, and $T_{0,n}(a,b)=T(a,b;p,n)$ if $p(x)=\log x$. The means $T_{q,n}$ belong to a large family of means called the \emph{hypergeometric means}, which have been studied by B. C. Carlson and others \cite{bc, c1, ct}. It is easy to see that $T_{q,n}$ is homogeneous. In particular, we have
$$T_{0,2}(a,b)=A(a,b),$$
$$T_{-2,2}(a,b)=G(a,b),$$
$$T_{2,2}(a,b)=Q(a,b).$$
Furthermore, the Toader means are related to the complete elliptic integrals: for $a\geq b>0$,
$$
T_{-1,2}(a,b)=\frac{\pi a}{2\ellipK(\sqrt{1-(b/a)^2})}
$$
and
$$
T_{1,2}(a,b)=\frac{2a}{\pi}\ellipE(\sqrt{1-(b/a)^2}).
$$

In 1997 S.-L. Qiu and J.-M. Shen \cite{qs2} proved that, for all $a,b>0$ with $a\neq b$,
$$
M_{3/2}(a,b)<T_{1,2}(a,b).
$$
This inequality had been conjectured by M. Vuorinen \cite{v}.  H. Alzer and S.-L. Qiu \cite{alq1} proved the following best possible  power mean upper bound:
$$
T_{1,2}(a,b)<M_{\log2/\log(\pi/2)}(a,b).
$$
Very recently, Y.-M. Chu and his collaborators \cite{cw,cw2,cwq} obtained several bounds for $T_{1,2}$ with respect to some combinations of various means.

\subsection{Seiffert means} The \emph{Seiffert means} $S_1$ and $S_2$ are defined by
$$
S_1(a,b)\equiv\frac{a-b}{2\arcsin\frac{a-b}{a+b}},\ a\neq b,\ \ S_1(a,a)=a,
$$
and
$$
S_2(a,b)\equiv\frac{a-b}{2\arctan\frac{a-b}{a+b}},\ a\neq b,\ \ S_2(a,a)=a.
$$

It is well known that
$$
\sqrt[3]{G^2A}<L<\frac{2G+A}{3}.
$$
S\'andor proved similar results for Seiffert means \cite{s2, s3}:
\begin{equation}\label{sandor1}
\sqrt[3]{A^2G}<S_1<\frac{G+2A}{3}<I
\end{equation}
and
\begin{equation}\label{sandor2}
\sqrt[3]{Q^2A}<S_2<\frac{A+2Q}{3}.
\end{equation}

The inequalities (\ref{sandor1}) and (\ref{sandor2}) are special cases of more general results obtained by E. Neuman and J. S\'andor \cite{ns4, ns5}.

\subsection{Extended means}
Let $a,b\in(0,\infty)$ be distinct and $s,t\in\R\setminus\{0\}$, $s\neq t$. We define the \emph{extended mean} \cite{st} with parameters $s$ and $t$ by
$$
E_{s,t}(a,b)\equiv\left(\frac{t}{s}\frac{a^s-b^s}{a^t-b^t}\right)^{1/(s-t)},
$$
and also
$$
E_{s,s}(a,b)\equiv \exp\left(\frac{1}{s}+\frac{a^s\log a-b^s\log b}{a^s-b^s}\right),
$$
$$
E_{s,0}(a,b)\equiv\left(\frac{a^s-b^s}{s\log(x/y)}\right)^{1/s}\quad {\rm and}\quad E_{0,0}(a,b)\equiv\sqrt{ab}.
$$
We see that all the classical means belong to the family of extended means. For example, $E_{2,1} = A$,  $E_{0,0} = G$, $E_{-1,-2}= H$,  and $E_{1,0}= L$, and, more generally, $M_\lambda = E_{2\lambda,\lambda}$ for $\lambda\in\R$. The reader is referred to the survey \cite{q} for many interesting results on the extended mean.

In 2002 P. A. H\"ast\"o \cite{ha} studied a certain monotonicity property of ratios of  extended means and Seiffert means, which he called a \emph{strong inequality}. These strong inequalities were shown to be related to the so-called \emph{relative metric} \cite{ha2, ha3}.

\subsection{Means and the circular and hyperbolic functions}
It is easy to check the following identities:
\begin{equation}\label{AGsin}
A(1+\sin x,1-\sin x)=1,\ \ G(1+\sin x,1-\sin x)=\cos x,
\end{equation}
\begin{equation}\label{QS1sin}
Q(1+\sin x,1-\sin x)=\sqrt{1+\sin^2x},\ \ S_1(1+\sin x,1-\sin x)=\frac{\sin x}{x},
\end{equation}
\begin{equation}\label{AGQcosh}
A(e^x,e^{-x})=\cosh x,\ \ G(e^x,e^{-x})=1, \ \ Q(e^x,e^{-x})=\sqrt{\cosh 2x},
\end{equation}
\begin{equation}\label{LIex}
L(e^x,e^{-x})=\frac{\sinh x}{x},\ \ I(e^x,e^{-x})=e^{x\coth x-1},
\end{equation}
\begin{equation}\label{S12ex}
S_1(e^x,e^{-x})=\frac{\sinh x}{\arcsin(\tanh x)},\ \ S_2(e^x,e^{-x})=\frac{\sinh x}{\arctan(\tanh x)}.
\end{equation}
One can get many inequalities by combining the above identities and inequalities between means. For example,
combining  (\ref{sandor1}) and (\ref{QS1sin})  we have
$$
\sqrt[3]{\cos x}<\frac{\sin x}{x}<\frac{\cos x+2}{3},
$$
where the second inequality is the well-known Cusa-Huygens inequality,
and
combining (\ref{sandor2}), (\ref{AGQcosh}), and (\ref{S12ex}), we have
$$
\sqrt[3]{(\cosh2x)(\cosh x)}<\frac{\sinh x}{\arctan(\tanh x)}<\frac{\cosh x+2\sqrt{\cosh 2x}}{3}.
$$
More inequalities on mean values and trigonometric and hyperbolic functions can be found in \cite{n9, s5, s4, y, zcl} and references therein.

\subsection{Means and hypergeometric functions}

In 2005 K. C. Richards \cite{ri} obtained sharp power mean bounds for the hypergeometric function:
Let $0<a,b\leq1$ and $c>\max\{-a,b\}$. If $c\geq\max\{1-2a,2b\}$, then
$$
M_\lambda(1-b/c;1,1-r)\leq F(-a,b;c;r)^{1/a}
$$
if and only if $\lambda\leq\frac{a+c}{1+c}$. If $c\leq\min\{1-2a,2b\}$, then
$$
M_\mu(1-b/c;1,1-r)\leq F(-a,b;c;r)^{1/a}
$$
if and only if $\mu\geq\frac{a+c}{1+c}$.
These inequalities generalize earlier results proved by Carlson \cite{c2}.

For hypergeometric functions of form $F(1/2-s,1/2+s;1;1-r^p)^q$, J. M. Borwein et al \cite{bbg} exhibited explicitly iterations similar to the arithmetic-geometric mean.  R. W. Barnard et al \cite{br} presented  sharp bounds for hypergeometric analogs of the arithmetic-geometric mean as follows: For $0<\alpha\leq1/2$ and $p>0$,
$$
M_\lambda(\alpha;1,r)\leq{F(\alpha,1-\alpha;1;1-r^p)^{-1/(\alpha p)}}\leq M_\mu(\alpha;1,r)
$$
if and only if $\lambda\leq0$ and $\mu\geq p(1-\alpha)/2$.

Some other inequalities involving hypergeometric functions and bivariate means can be found in the very recent survey \cite{brt}.

For any two power means $M_\lambda$ and $M_\mu$, a function $f$ is called $M_{\lambda,\mu}$-\emph{convex} if it satisfies
$$
f(M_\lambda(x,y))\leq M_\mu(f(x),f(y)).
$$
Recently many authors have proved that the zero-balanced Gaussian hypergeometric function is $M_{\lambda,\lambda}$-convex  when $\lambda\in\{-1, 0, 1\}$. For details see
\cite{avv8}, \cite{bpv}, \cite{bar4}, and \cite{cwzq}. \'A. Baricz \cite{bar2} generalized these results to the $M_{\lambda,\lambda}$-convexity of zero-balanced Gaussian hypergeometric functions with respect to a power mean for $\lambda\in[0,1]$.  X.-H. Zhang et al \cite{zwc} extended these results to the case of $M_{\lambda,\mu}$-convexity with respect to two power means: For all $a,b>0$, $\lambda\in(-\infty,1]$, and $\mu\in[0,\infty)$ the hypergeometric function
$F(a,b;a+b;r)$ is $M_{\lambda,\mu}$-convex on $(0,1)$.

The following interesting open problem is presented by \'A. Baricz \cite{bar8}:

\medskip

\noindent{\bf Open problem.} {\it If $m_1$ and $m_2$ are bivariate means, then find conditions on $a_1,a_2 > 0$ and $c> 0$ for
which the inequality
$$
m_1(F_{a_1}(r), F_{a_2}(r))\leq(\geq)F_{m_2(a_1,a_2)}(r)
$$
holds true for all $ r \in(0, 1)$, where $F_a(r)=F(a,c-a;c;r)$.}

\subsection{Means and quasiconformal analysis}

Special functions have always played an important role in the distortion theory of quasiconformal mappings.
G. D. Anderson, M. K. Vamanamurthy, and M. Vuorinen \cite{avv2} have systematically investigated classical special functions and their extensive applications in the theory of conformal invariants and quasiconformal mappings. Some functional inequalities
for special functions in quasiconformal mapping theory involve the arithmetic mean, geometric mean, or harmonic mean.
For example, for the well-known Gr\"otzsch ring function $\mu$ and the Hersch-Pfluger distortion function $\varphi_K$, the following inequalities hold for all $s,t\in(0,1)$ with $s\neq t$:
$$
\sqrt{\mu(s)\mu(t)}<\mu(\sqrt{st}),
$$
and
$$
\sqrt{\varphi_K(s)\varphi_K(t)}<\varphi_K(\sqrt{st})  \ \  \mbox{for}\  K>1.
$$
Recently, G.-D. Wang, X.-H. Zhang, and Y.-M. Chu \cite{wzc3, wzc2} have extended these inequalities as follows:
$$
M_\lambda(\mu(s),\mu(t))<\mu(M_\lambda(s,t))\ \ \mbox{if and only if }\ \lambda\leq0,
$$
$$
M_\lambda(\varphi_K(s),\varphi_K(t))<\varphi_K(M_\lambda(s,t))\ \ \mbox{if and only if }\ \lambda\geq0 \ {\rm and }\ K>1,
$$
and
$$
M_\lambda(\varphi_K(s),\varphi_K(t))>\varphi_K(M_\lambda(s,t))\ \ \mbox{if and only if }\ \lambda\geq0 \ {\rm and }\ 0<K<1.
$$

Some similar results for the generalized Gr\"otzsch function, generalized modular function, and other special functions related to quasiconformal analysis can be found in \cite{qqwc, wqzc, wzc4, wzj1, wzj2}.

\section{Epilogue and a view toward the future}

In earlier work we have listed many open problems.  See especially \cite[ pp. 128--131]{avv6} and \cite[p. 478]{avv2}. Many of these problems are still open. In Sections 4, 6, and 7 above we have also mentioned some open problems.

 Finally, we wish to suggest some ideas for further research.  In a frequently cited paper \cite{lt}  P. Lindqvist introduced in 1995 the notion of \emph{generalized trigonometric functions} such as $\sin_p$, and presently there is a large body of literature about this topic.  For the case $p=2$ the classical functions are obtained.  In 2010 R. J. Biezuner et al \cite{bem} developed a practical numerical method for computing values of $\sin_p$.  Recently, S. Takeuchi \cite{ta} has gone a step further, introducing functions depending on two parameters $p$ and $q$ that reduce to the $p$-functions of Lindqvist when $p=q$.  In \cite{bv2, bv3, kvz} the authors have continued the study of this family of  generalized  functions, and have suggested that many properties of classical functions have a counterpart in this more general setting.   It would be natural to generalize the properties of trigonometric functions cited in this survey to the $(p,q)$-trigonometric functions of Takeuchi.

\begin{acknowledgement}
The authors wish to thank \'A. Baricz, C. Berg, E. A. Karatsuba, C. Mortici, E. Neuman, H. L. Pedersen, S. Ponnusamy, and G. Tee  for careful reading of this paper and for many corrections and suggestions. The research of Matti Vuorinen was supported by the Academy of Finland, Project 2600066611. Xiaohui Zhang is indebted to the Finnish National Graduate School of Mathematics and its Applications for financial support.
\end{acknowledgement}

%
%

\begin{thebibliography}{999.}


\bibitem{as} Abramowitz, M., Stegun, I.A.: Handbook of Mathematical Functions with Formulas, Graphs, and Mathematical Tables. Dover Publ., New York (1970)

\bibitem{aa} Adell, J.A., Alzer, H.: A monotonicity property of Euler's gamma function. Publ. Math. Debrecen \textbf{78}, 443--448  (2011)

\bibitem{alzer}Alzer, H.: On some inequalities for the gamma and psi functions. Math. Comput. \textbf{66}, 373--389 (1997)

\bibitem{a1}Alzer, H.: Inequalities for the gamma function.  Proc. Amer. Math. Soc. \textbf{128}, 141--147  (1999)

\bibitem{a3}Alzer, H.: Inequalities for the volume of the unit ball in $\mathbb{R}^n$.  J. Math. Anal. Appl. \textbf{252}, 353--363  (2000)

\bibitem{a0}Alzer, H.: Sharp inequalities for digamma and polygamma functions. Forum Math. \textbf{16}, 181--221  (2004)

\bibitem {a4}Alzer, H.: Inequalities for the volume of the unit ball in $\mathbb{R}^n$ II. Mediterr. J. Math. \textbf{5}, 395--413 (2008)

\bibitem{a}Alzer, H.: Inequalities for the harmonic numbers. Math. Z. \textbf{267}, 367--384  (2011)

\bibitem{alb}Alzer, H., Batir, N.: Monotonicity properties of the gamma function.  Appl. Math. Lett. \textbf{20}, 778--781 (2007)

\bibitem{alq1}Alzer, H., Qiu, S.-L.: Inequalities for means in two variables.  Arch. Math. (Basel) \textbf{80}, 201--215  (2003)

\bibitem{alq}Alzer, H., Qiu, S.-L.: Monotonicity theorems and inequalities for the complete elliptic integrals.  J. Comput. Appl. Math. \textbf{172}, 289--312  (2004)

\bibitem{aq} Anderson, G.D., Qiu, S.-L.: A monotoneity property of the gamma function.  Proc. Amer. Math. Soc. \textbf{125}, 3355--3362 (1997)

\bibitem{aqvv} Anderson, G.D., Qiu, S.-L., Vamanamurthy, M.K., Vuorinen, M.: Generalized
   elliptic integrals and modular equations.  Pacific J. Math. \textbf{192}, 1--37 (2000)

\bibitem{avv6} Anderson, G.D., Vamanamurthy, M.K., Vuorinen, M.: Special functions of quasiconformal theory.  Expo. Math. \textbf{7}, 97--136  (1989)

\bibitem{avv7} Anderson, G.D., Vamanamurthy, M.K., Vuorinen, M.: Functional inequalities for hypergeometric functions and complete elliptic integrals.  SIAM J. Math. Anal. \textbf{23}, 512--524 (1992)

\bibitem{avv9} Anderson, G.D., Vamanamurthy, M.K., Vuorinen, M.: Hypergeometric
    functions and elliptic integrals. In:  Srivastava, H.M., Owa, S. (eds.) Current Topics in Analytic Function Theory, pp. 48-85. World Scientific Publ. Co.  (1992)

\bibitem {avv1} Anderson, G.D., Vamanamurthy, M.K., Vuorinen, M.: Inequalities for quasiconformal mappings in space. Pacific J. Math. \textbf{160}, 1--18  (1993)

\bibitem {avv2} Anderson, G.D., Vamanamurthy, M.K., Vuorinen, M.: Conformal Invariants, Inequalities, and Quasiconformal Maps. J. Wiley, New York (1997)

\bibitem {avv3} Anderson, G.D., Vamanamurthy, M.K., Vuorinen, M.: Topics in special functions. In: Papers on Analysis: A volume dedicated to Olli Martio on the occasion of his 60th birthday, Report Univ. Jyv\"askyl\"a \textbf{83}, 5--26 (2001)

\bibitem {avv4} Anderson, G.D., Vamanamurthy, M.K., Vuorinen, M.: Monotonicity rules in calculus. Amer. Math. Monthly \textbf{133}, 805--816  (2006)

\bibitem {avv5} Anderson, G.D., Vamanamurthy, M.K., Vuorinen, M.: Topics in special functions  II.   Conform. Geom. Dyn. \textbf{11}, 250--271  (2007)

\bibitem {avv8} Anderson, G.D., Vamanamurthy, M.K., Vuorinen, M.: Generalized convexity and inequalities. J. Math. Anal. Appl. \textbf{335}, 1294--1308 (2007)


\bibitem {av1}Anderson, G.D., Vuorinen, M.: Reflections on Ramanujan's mathematical gems.
    Math. Newsl. \textbf{19}, 87--108 (2010) Available via arXiv:1006.5092v1 [math.CV]

\bibitem {ab} Andr\'{a}s, S., Baricz, \'{A}.: Bounds for complete elliptic integrals of the first kind.  Expo. Math. \textbf{28}, 357--364 (2010)

\bibitem {aar}Andrews, G., Askey,  R., Roy, R.: Special Functions. Encyclopedia of Mathematics
and its Applications  Vol.\,71. Cambridge Univ. Press  (1999)

\bibitem {bnpv}  Balasubramanian, R., Naik,  S.,  Ponnusamy, S.,  Vuorinen, M.: Elliott's identity and hypergeometric functions.   J. Math. Anal. Appl. \textbf{271}, 232--256 (2002)

\bibitem {bpv}  Balasubramanian, R.,   Ponnusamy, S.,  Vuorinen, M.: Functional inequalities for the quotients of hypergeometric functions.  J. Math. Anal. Appl. \textbf{218}, 256--268 (1998)

\bibitem  {bapi} Barbu, C., Pi\c{s}coran,  L.-I.: On Panaitopol and Jordan type inequalities.  Unpublished manuscript

\bibitem {bar3}   Baricz, \'{A}.: Landen-type inequalities for Bessel functions.  Comput. Methods Funct. Theory \textbf{5}, 373-379 (2005)

\bibitem {bar5} Baricz, \'{A}.: Functional inequalities involving special functions. J. Math. Anal. Appl. \textbf{319}, 450--459 (2006)

\bibitem {bar7}  Baricz, \'{A}.: Tur\'an type inequalities for generalized complete elliptic integrals.  Math. Z. \textbf{256}, 895--911  (2007)

\bibitem {bar9} Baricz, \'{A}.: Functional inequalities involving special functions II.  J. Math. Anal. Appl. \textbf{327}, 1202--1213  (2007)

\bibitem {bar2}  Baricz, \'{A}.: Convexity of the zero-balanced Gaussian hypergeometric functions with respect to H\"older means.  JIPAM. J. Inequal. Pure Appl. Math.  \textbf{8}, Article 40, 9 pp.  (2007)

\bibitem  {bar1}  Baricz, \'{A}.: Jordan-type inequalities for generalized Bessel functions.  JIPAM. J. Inequal. Pure Appl. Math. \textbf{9}, Article 39, 6 pp. (2008)

\bibitem {bar6}  Baricz, \'{A}.: Functional inequalities involving Bessel and modified Bessel functions of the first kind.   Expo. Math. \textbf{26},  279--293 (2008)

\bibitem {bar8}  Baricz, \'{A}.: Tur\'an type inequalities for hypergeometric functions. Proc. Amer. Math. Soc. \textbf{136}, 3223--3229 (2008)

\bibitem {bar4} Baricz, \'{A}.: Generalized Bessel functions of the first kind. Lecture Notes in Mathematics 1994, Springer-Verlag, Berlin  (2010)
    
\bibitem {bar10} Baricz, \'{A}.: Landen inequalities for special functions.  Proc. Amer. Math. Soc. (to appear). Available via arXiv:1301.5255 [math.CA]

\bibitem {bs}  Baricz, \'{A}., S\'{a}ndor, J.: Extensions of the generalized Wilker inequality to Bessel functions.  J. Math. Inequal. \textbf{2}, 397--406 (2008)

\bibitem {bvv}  Baricz, \'{A}., Vesti, J., Vuorinen, M.: On Kaluza's sign criterion for reciprocal power series.  Ann. Univ. Mariae Curie-Sk{\l}odowska Sect A \textbf{65}, 1--16 (2011)

\bibitem {bw}  Baricz, \'{A}., Wu, S.-H.: Sharp Jordan type inequalities for Bessel functions.  Publ. Math. Debrecen \textbf{74}, 107--126 (2009)

\bibitem {bw2}  Baricz, \'{A}., Wu, S.-H.: Sharp exponential Redheffer-type inequalities for Bessel functions.  Publ. Math. Debrecen \textbf{74}, 257--278 (2009)


\bibitem {br} Barnard, R.W., Richards, K.C.: On inequalities for hypergeometric analogues of the arithmetic-geometric mean.  JIPAM. J. Inequal. Pure Appl. Math.  \textbf{8}, Article 65, 5 pp.  (2007)

\bibitem {brt} Barnard, R.W., Richards, K.C., Tiedeman, H.C.: A survey of some bounds for Gauss' hypergeometric function and related bivariate means.  J. Math. Inequal. \textbf{4}, 45--52 (2010)

\bibitem {bat3}Batir, N.: On some properties of digamma and polygamma functions.  J. Math. Anal. Appl. \textbf{328}, 452--465  (2007)

\bibitem {bat2}Batir, N.: On some properties of the gamma function.  Expo. Math. \textbf{26}, 187--196 (2008)

\bibitem {bat}Batir, N.: Sharp inequalities for factorial $n$.  Proyecciones \textbf{27}, 97--102 (2008)

\bibitem {bat4}Batir, N.: Improving Stirling's formula.  Math. Commun. \textbf{16}, 105--114 (2011)

\bibitem {bes} Becker, M., Stark, E.L.: On a hierarchy of quolynomial inequalities for $\tan x.$   Univ. Beograd. Publ. Elektrotehn. Fak. Ser. Mat. Fiz.  No. 602--633 (1978), 133--138 (1979)

\bibitem {bp1}  Berg, C., Pedersen, H.L.: A completely monotone function related to the
    gamma function. J. Comput. Appl. Math. \textbf{133}, 219--230  (2001)

\bibitem {bp2}  Berg, C., Pedersen, H.L.: Pick functions related to the gamma function.  Rocky Mountain J. Math.  \textbf{32}, 507--525  (2002)

\bibitem {bp3}  Berg, C., Pedersen, H.L.: A one-parameter family of Pick functions defined by the gamma function and related to the volume of the unit ball in n-space.  Proc. Amer. Math. Soc. \textbf{139}, 2121--2132 (2011)

\bibitem {bp}  Berg, C., Pedersen, H.L.: A completely monotonic function used in an inequality of Alzer.  Comput. Methods Funct. Theory \textbf{12}, 329--341 (2012)
    
\bibitem {bm} Berinde, V., Mortici, C.: New sharp estimates of the generalized Euler-Mascheronoi constant. Math. Inequal. Appl. \textbf{16}, 279--288 (2013)

\bibitem {bt}  Berndt, B.C.: Ramanujan's Notebooks, Part II. Springer-Verlag, New York (1987)

\bibitem {bv2} Bhayo, B.A., Vuorinen, M.: On generalized trigonometric functions with two parameters.   J. Approx. Theory \textbf{164}, 1415--1426  (2012)

\bibitem {bv} Bhayo, B.A., Vuorinen, M.: On generalized complete elliptic integrals and modular functions.  Proc. Edinburgh Math. Soc. \textbf{55}, 591--611 (2012)

\bibitem {bv3} Bhayo, B.A., Vuorinen, M.: Inequalities for eigenfunctions of the $p$-Laplacian.  Available via
   	arXiv:1101.3911v3 [math.CA]

\bibitem {bk} Biernacki, M.,  Krzy\.z, J.: On the monotonicity of certain functionals in the theory of analytic functions.  Ann. Univ. M. Curie-Sk{\l}odowska \textbf{2}, 134--145 (1995)

\bibitem {bem} Biezuner, R.J., Ercole, G.,  Martins, E.M.: Computing the $\sin_p$ function via the inverse power method.  Comput. Methods Appl. Math. \textbf{2}, 129--140 (2012)

\bibitem {bh} B\"ohm,  J., Hertel, E.: Polyedergeometrie in n-Dimensionalen R\"aumen Konstanter Kr\"ummung.  Birkh\"auser, Basel  (1981)

\bibitem {bbg} Borwein, J.M.,  Borwein, P.B., Garvan, F.: Hypergeometric analogues of the arithmetic-geometric mean iteration.  Constr. Approx. \textbf{9}, 509--523 (1993)

\bibitem {bc} Brenner, J.L., Carlson,  B.C.: Homogeneous mean values: weights and asymptotics.  J. Math. Anal. Appl. \textbf{123}, 265--280 (1987)

\bibitem {bu} Burnside,  W.: A rapidly convergent series for $\log N!$.  Messenger Math. \textbf{46}, 157--159 (1917)

\bibitem  {bf} Byrd, P.F.,  Friedman, M.D.:
    Handbook of Elliptic Integrals for Engineers and Scientists. 2nd ed., Die Grundlehren der Math. Wiss. 67, Springer-Verlag, Berlin  (1971)

\bibitem {c1}  Carlson, B.C.: A hypergeometric mean value.  Proc. Amer. Math. Soc. \textbf{16}, 759--766 (1965)

\bibitem {c2}   Carlson, B.C.: Some Inequalities for hypergeometric functions.  Proc. Amer. Math. Soc. \textbf{16}, 32--39 (1966)

\bibitem {c}   Carlson, B.C.: Inequalities for a symmetric elliptic integral.  Proc. Amer. Math. Soc. \textbf{25}, 698--703 (1970)

\bibitem {ct}  Carlson, B.C., Tobey, M.D.:  A property of the hypergeometric mean value.  Proc. Amer. Math. Soc. \textbf{19}, 255--262 (1968)

\bibitem {cn2} Chen, C.-P.: Inequalities for the Euler-Mascheroni constant.  Appl. Math. Lett. \textbf{23}, 161--164 (2010)

\bibitem {cn}  Chen, C.-P.: Sharpness of Negoi's inequality for the Euler-Mascheroni constant.  Bull. Math. Anal. Appl. \textbf{3}, 134--141 (2011)

\bibitem {cc}  Chen, C.-P., Cheung, W.-S.: Sharp Cusa and Becker-Stark inequalities.  J. Inequal. Appl. \textbf{2011}, Article 136, 6 pp. (2011)

\bibitem {cc2} Chen, C.-P., Cheung, W.-S.: Sharpness of Wilker and Huygens type inequalities.  J. Inequal. Appl. \textbf{2012}, Article 72, 11 pp. (2012)

\bibitem {ccw} Chen, C.-P., Cheung, W.-S., Wang, W.-S.: On Shafer and Carlson inequalities.  J. Inequal. Appl. \textbf{2011}, Article ID 840206, 10 pp.  (2011)

\bibitem {cd} Chen, C.-P., Debnath, L.: Sharpness and generalization of Jordan's inequality and its application.  Appl. Math. Lett. \textbf{25}, 594--599 (2012)

\bibitem {cm} Chen, C.-P., Mortici, C.: Generalization and sharpness of Carlson's inequality for the inverse cosine function.  Unpublished manuscript

\bibitem {czq} Chen, C.-P., Zhao, J.-W.,  Qi, F.: Three inequalities involving hyperbolically trigonometric functions.   RGMIA Res. Rep. Coll. \textbf{6}(3), Article 4, 437-443 (2003)

\bibitem {chl} Chlebus, E.: A recursive scheme for improving the original rate of convergence to the Euler-Mascheroni constant.  Amer. Math. Monthly \textbf{118}, 268--274 (2011)

\bibitem {cwzq} Chu, Y.-M.,   Wang, G.-D., Zhang, X.-H.,  Qiu, S.-L.: Generalized convexity and inequalities involving special functions.  J. Math. Anal. Appl. \textbf{336}, 768-776 (2007)

\bibitem {cw}  Chu, Y.-M., Wang, M.-K.: Optimal Lehmer mean bounds for the Toader mean.  Results Math. \textbf{61}, 223--229 (2012)

\bibitem {cw2} Chu, Y.-M., Wang, M.-K.: Inequalities between arithmetic-geometric, Gini, and Toader Means.  Abstr. Appl. Anal. \textbf{2012}, Article ID 830585, 11 pp.  (2012)

\bibitem {cwq} Chu, Y.-M., Wang, M.-K., Qiu, S.-L.: Optimal combinations bounds of root-square and arithmetic means for Toader mean.  Proc. Indian Acad. Sci. Math. Sci. \textbf{122}, 41--51 (2012)

\bibitem {cwq2} Chu, Y.-M., Wang, M.-K., Qiu, Y.-F.: On Alzer and Qiu's conjecture for complete elliptic integral and inverse hyperbolic tangent function.  Abstr. Appl. Anal.  \textbf{2011}, Article ID 697547, 7 pp. (2011)

\bibitem {cwqj} Chu, Y.-M., Wang, M.-K., Qiu, S.-L., Jiang, Y.-P.: Bounds for complete elliptic integrals of the second kind with applications.  Comput. Math. Appl. \textbf{63}, 1177--1184  (2012)

\bibitem {d} DeTemple, D.W.: Convergence to Euler's constant.  Amer. Math. Monthly \textbf{100}, 468--470 (1993)

\bibitem {el} Elbert,\'A.,  Laforgia, A.: On some properties of the gamma function.    Proc. Amer. Math. Soc. \textbf{128}, 2667--2673  (2000)

\bibitem {e} Elliott, E.B.: A formula including Legendre's  $\ellipE\ellipK'+\ellipK\ellipE'-\ellipK\ellipK'=\frac{1}{2}\pi$.  Messenger of Math. \textbf{33}, 31--40 (1904)

\bibitem {f}  Fink, A.M.: Two inequalities.  Univ. Beograd. Publ. Elektrotehn. Fak. Ser. Mat. \textbf{6}, 49--50 (1995)

\bibitem {g} Ge, H.-F.: New sharp bounds for the Bernoulli numbers and refinement of Becker-Stark inequalities.  J. Appl. Math. \textbf{2012}, Article ID 137507, 7 pp.  (2012)

\bibitem {gcq} Guo,  B.-N., Chen,  R.-J.,  Qi, F.: A class of completely monotonic functions involving the polygamma functions.   J. Math. Anal. Approx. Theory \textbf{1}, 124--134 (2006)

\bibitem  {gq} Guo,  B.-N.,  Qi, F.: Some bounds for the complete elliptic integrals of the first and second kinds. Math. Inequal. Appl. \textbf{14}, 323--334 (2011)

\bibitem {ha} H\"ast\"o, P.A.: A monotonicity property of ratios of symmetric homogeneous means.  JIPAM. J. Inequal. Pure Appl. Math. \textbf{3}, Article 71, 23 pp.  (2002)

\bibitem {ha2} H\"ast\"o, P.A.: A new weighted metric: the relative metric I.  J. Math. Anal. Appl.  \textbf{274}, 38--58 (2002)

\bibitem {ha3}H\"ast\"o, P.A.: A new weighted metric: the relative metric II.  J. Math. Anal. Appl.  \textbf{301}, 336--353 (2005)

\bibitem {hlvv} Heikkala, V., Lind\'en, H.,  Vamanamurthy, M. K., Vuorinen, M.: Generalized
elliptic integrals and the Legendre M -function.  J. Math. Anal. Appl. \textbf{338}, 223--243 (2008)

\bibitem  {hvv} Heikkala, V., Vamanamurthy, M. K., Vuorinen, M.: Generalized elliptic
integrals.  Comput. Methods Funct. Theory \textbf{9}, 75--109 (2009)

\bibitem {hu} Hua, Y.: Refinements and sharpness of some new Huygens type inequalities.  J. Math. Inequal. \textbf{6},  493--500 (2012)

\bibitem{hncq}Huo, Z.-H., Niu,  D.-W., Cao,  J.,  Qi, F.: A generalization of Jordan's inequality and an application.  Hacet. J. Math. Stat. \textbf{40}, 53--61 (2011)

\bibitem {hs} Huygens, C.: Oeuvres Completes 1888-1940.   Soci\'et\'e Hollondaise des Science, Haga

\bibitem {i} Iv\'ady, P.: A note on a gamma function inequality.  J. Math. Inequal. \textbf{3}, 227--236 (2009)

\bibitem {id} Ibrahim, A., Dragomir, S.\,S.: Power series inequalities via Buzano's result and applications. Integral Transforms Spec. Funct.  \textbf{22}, 867--878  (2011)

\bibitem {idcd}Ibrahim, A., Dragomir, S.\,S., Cerone, P., Darus, M.: Inequalities for power series with positive coefficients. J. Inequal. Spec. Funct. \textbf{3}, 1--15 (2012)

\bibitem {idd}Ibrahim, A., Dragomir, Darus, M.: Some inequalities for power series with applications. Integral Transforms Spec. Funct. iFirst, 1-13 (2012)

\bibitem {kk1} Kalmykov, S.I., Karp, D.B.: Log-concavity for series in reciprocal gamma functions and applications. Integral Transforms Spec. Funct. (2013). Available via arXiv:1206.4814v1 [math.CA]

\bibitem {kk2} Kalmykov, S.I., Karp, D.B.: Log-convexity and log-concavity for series in gamma ratios and applications. Available via  arXiv:1211.2882v2 [math.CA]

\bibitem {kal} Kaluza, T.: \"Uber die Koeffizienten reziproker Potenzreihen.  Math. Z. \textbf{28}, 161--170 (1928)

\bibitem {ka} Karatsuba, E.A.: On the asymptotic representation of the Euler gamma function by Ramanujan.   J. Comput. Appl. Math. \textbf{135}, 225--240 (2001)

\bibitem {kv}  Karatsuba, E.A.,  Vuorinen, M.: On hypergeometric functions and generalizations
of Legendre's relation.  J. Math. Anal. Appl. \textbf{260}, 623--640 (2001)

\bibitem {ks} Karp, D., Sitnik, S.M.: Inequalities and monotonicity of ratios for generalized
    hypergeometric function.  J. Approx. Theory \textbf{161}, 337--352 (2009)

\bibitem {kmsv} Kl\'{e}n, R., Manojlovic,  V., Simi\'{c}, S., Vuorinen, M.: Bernoulli inequality and hypergeometric functions.  Proc. Amer. Math. Soc. (to appear)

\bibitem {kmv} Kl\'en,  R., Manojlovi\'c, V., Vuorinen, M.: Distortion of normalized quasiconformal mappings. Available via arXiv:0808.1219 [math.CV]

\bibitem  {kvv} Kl\'en,  R., Visuri,  M., Vuorinen, M.: On Jordan type inequalities for hyperbolic functions.  J. Inequal. Appl. \textbf{2010}, Article ID 362548, 14 pp.  (2010)
    
\bibitem  {kvz} Kl\'en,  R., Vuorinen, M., Zhang, X.-H.: Inequalities for the generalized trigonometric and hyperbolic functions. Available via arXiv:1210.6749 [math.CA]

\bibitem  {k1} Kouba, O.: New bounds for the identric mean of two arguments.  JIPAM. J. Inequal. Pure Appl. Math. \textbf{9}, Article 71.  (2008)

\bibitem  {k2} Kouba, O.: Bounds for the ratios of differences of power means in two arguments. Math. Inequal. Appl. (to appear). Available via arXiv:1006.1460v1 [math.CA]

\bibitem {kp} Koumandos, S., Pedersen, H.L.: On the asymptotic expansion of the logarithm of Barnes triple gamma function.  Math. Scand. \textbf{105}, 287--306 (2009)

\bibitem {ko}Kuo, M.-K.: Refinements of Jordan's inequality.  J. Inequal. Appl. \textbf{2011}:130, 6 pp.  (2011)

\bibitem {l} Lazarevi\'c, I.: Neke nejednakosti sa hiperbolickim funk\u{c}ijama.  Univ. Beograd. Publ. Elektrotehn. Fak. Ser. Mat. \textbf{170}, 41--48 (1966)

\bibitem  {lv} Lehto, O.,  Virtanen, K.I.: Quasiconformal Mappings in the Plane, 2nd ed., Grundlehren Math. Wiss., Band 126, Springer-Verlag, New York (1973)

\bibitem  {lt} Lindqvist, P.: Some remarkable sine and cosine functions.  Ric. Mat. \textbf{44}, 269--290 (1995)

\bibitem {ll} Li, J.-L., Li, Y.-L.: On  the strengthened Jordan's inequality.  J. Inequal. Appl.  \textbf{2007}, Article ID 74328, 8 pp.  (2007)

\bibitem  {lwc} Lv, Y.-P., Wang,  G.-D., Chu, Y.-M.: A note on Jordan type inequalities for hyperbolic functions. Appl. Math. Lett. \textbf{25}, 505--508 (2012)

\bibitem {mqzc} Ma, X.-Y., Qiu,  S.-L., Zhong, G.-H., Chu, Y.-M.: Some inequalities for the generalized linear distortion function.  Appl. Math. J. Chinese Univ. Ser. B \textbf{27}, 87--93 (2012)

\bibitem {mqzc2} Ma, X.-Y., Qiu,  S.-L., Zhong, G.-H., Chu, Y.-M.: The H\"{o}lder continuity and submultiplicative properties of the modular function. Appl. Math. J. Chinese Univ. Ser. A \textbf{27}, 481--487 (2012)

\bibitem  {maa}Mahmoud, M., Alghamdi,  M.A., Agarwal, R.P.: New upper bounds of $n!$.  J. Inequal. Appl. \textbf{2012}, Article 27, 9 pp. (2012)

\bibitem  {ms} Miller, K.S.,  Samko, S.G.: Completely monotonic functions.   Integral Transforms Spec. Funct. \textbf{12}, 389--402 (2001)

\bibitem  {mi} Mitrinovi\'c, D.S.: Analytic inequalities. Springer-Verlag, Berlin (1970)

\bibitem {mo} Mori, A.: On an absolute constant in the theory of quasiconformal mappings.  J. Math. Soc. Japan \textbf{8}, 156--166 (1956)

\bibitem  {m2} Mortici, C.: Monotonicity properties of the volume of the unit ball in $\mathbb{R}^n$.  Optim. Lett. \textbf{4}, 457--464 (2010)

\bibitem  {m3} Mortici, C.: Very accurate estimates of the polygamma functions. Asymptot. Anal. \textbf{68}, 125--134 (2010)

\bibitem  {m4} Mortici, C.: On new sequences converging towards the Euler-Mascheroni constant.  Comput. Math. Appl. \textbf{59}, 2610--2614 (2010)

\bibitem {m5} Mortici, C.: New approximations of the gamma function in terms of the digamma function.  Appl. Math. Lett. \textbf{23}, 97--100 (2010)

\bibitem  {m} Mortici, C.: Ramanujan's estimate for the gamma function via monotonicity arguments.  Ramanujan J. \textbf{25}, 149--154 (2011)

\bibitem  {m8} Mortici, C.: Gamma function by $x^{x-1}$.   Carpathian J. Math. (to appear)

\bibitem {ni} Negoi, T.: A faster convergence to Euler's constant.  Math. Gaz. \textbf{83}, 487--489 (1999)

\bibitem {n2} Neuman, E.: Inequalities and bounds for generalized complete elliptic integrals.  J. Math. Anal. Appl. \textbf{373}, 203--213 (2011)

\bibitem {n4} Neuman, E.: On Wilker and Huygens type inequalites.  Math. Inequal. Appl. \textbf{15}, 271--279 (2012)

\bibitem {n}  Neuman, E.: Inequalities involving hyperbolic functions and trigonometric functions.  Bull. Int. Math. Virt. Instit. \textbf{2}, 87--92 (2012)

\bibitem {n9} Neuman, E.: A note on a certain bivariate mean.  J. Math. Inequal. \textbf{6}, 637--643 (2012)

\bibitem {ns4} Neuman, E., S\'{a}ndor, J.: On the Schwab-Borchardt mean. Math. Pannon. \textbf{14}, 253--266 (2003)

\bibitem {ns5} Neuman, E., S\'{a}ndor, J.: On the Schwab-Borchardt mean II. Math. Pannon. \textbf{17}, 49--59 (2006)

\bibitem {ns1} Neuman, E., S\'{a}ndor, J.: On some inequalities involving trigonometric and hyperbolic functions with emphasis on the Cusa-Huygens, Wilker, and Huygens inequalities.  Math. Inequal. Appl. \textbf{13}, 715--723 (2010)

\bibitem {ns}Neuman, E., S\'{a}ndor, J.: Optimal inequalities for hyperbolic and trigonometric functions.  Bull. Math. Anal. Appl. \textbf{3}, 177--181 (2011)

\bibitem {ncq}Niu,  D.-W., Cao,  J., Qi, F.: Generalizations of Jordan's inequality and concerned relations.  Politehn. Univ. Bucharest Sci. Bull. Ser. A Appl. Math. Phys. \textbf{72}, 85--98 (2010)

\bibitem {nhcq} Niu, D.-W., Huo, Z.-H., Cao, J., Qi, F.: A general refinement
    of Jordan's inequality and a refinement of L. Yang's inequality. Integral Transforms Spec. Funct. \textbf{19}, 157--164 (2008)

\bibitem {pz}Pan, W.-H.,  Zhu, L.: Generalizations of Shafer-Fink-type inequalities for the arc sine function.  J. Inequal. Appl. \textbf{2009}, Article ID 705317, 6 pp.  (2009)

\bibitem  {pi2} Pinelis, I.: L'Hospital rules for monotonicity and the Wilker-Anglesio inequality. Amer. Math. Monthly \textbf{111}, 905--909 (2004)

\bibitem {pv} Ponnusamy, S.,  Vuorinen, M.: Asymptotic expansions and inequalities for hypergeometric functions. Mathematika \textbf{44}, 278--301 (1997)

\bibitem {q} Qi, F.: The extended mean values: definition, properties, monotonicities, comparison, convexities, generalizations, and applications.  RGMIA Res. Rep. Coll. \textbf{5}, 19 pp. (2001)

\bibitem {qg} Qi, F., Guo, B.-N.: Monotonicity and logarithmic convexity relating to the volume of the unit ball.   Optim. Lett. (2012). Available via arXiv:0902.2509v1 [math.CA]

\bibitem {qng} Qi, F., Niu, D.-W., Guo, B.-N.: Refinements, generalizations, and applications of Jordan's inequality and related problems. J. Inequal. Appl. \textbf{2009}, Article ID 271923, 52 pp. (2009)

\bibitem {qqwc}Qiu, S.-L., Qiu,  Y.-F., Wang, M.-K., Chu, Y.-M.: H\"older mean inequalities for the generalized Gr\"otzsch ring and Hersch-Pfluger functions. Math. Inequal. Appl. \textbf{15}, 237--245 (2012)

\bibitem{qs2} Qiu,  S.-L., Shen, J.-M.: On two problems concerning means. J. Hangzhou Inst. Electronic Engg. \textbf{17}, 1-7 (1997)

\bibitem {qv2} Qiu,  S.-L., Vuorinen, M.: Landen inequalities for hypergeometric functions.  Nagoya Math. J. \textbf{154}, 31--56(1999)

\bibitem {qv3} Qiu, S.-L., Vuorinen, M.: Duplication inequalities for the ratios of hypergeometric functions.  Forum Math. \textbf{12}, 109--133 (2000)

\bibitem {qv} Qiu, S.-L., Vuorinen, M.: Some properties of the gamma and psi functions, with applications.  Math Comput. \textbf{74}, 723--742 (2004)

\bibitem {ra} Ramanujan, S.:  The Lost Notebook and Other Unpublished Papers, with an Introduction by George E. Andrews.  Naros Publishing House, New Delhi, Madras, Bombay (1988)

\bibitem {r}Redheffer, R.: Problem 5642. Amer. Math. Monthly \textbf{76}, 422 (1969)

\bibitem  {ri} Richards, K.C.: Sharp power mean bounds for the Gaussian hypergeometric function.  J. Math. Anal. Appl. \textbf{308}, 303--313 (2005)

\bibitem {s7}S\'{a}ndor, J.: Sur la fonction gamma. Publ. Centre Rech. Math. Pures (I) \textbf{21}, 4--7 (1989)

\bibitem {s2}S\'{a}ndor, J.: On certain inequalities for means III.  Arch. Math. (Basel) \textbf{76}, 34--40 (2001)

\bibitem {s3}S\'{a}ndor, J.: \"Uber zwei Mittel von Seiffert. Wurzel \textbf{36}, 104--107 (2002)

\bibitem {s1}S\'{a}ndor, J.: Two sharp inequalities for trigonometric and hyperbolic functions.  Math. Inequal. Appl. \textbf{15}, 409--413 (2012)

\bibitem {s5}S\'{a}ndor, J.: On Huygens' inequalities and the theory of means.  Int. J. Math. Math. Sci.  \textbf{2012}, Article ID 597490, 9 pp.  (2012)

\bibitem {s6}S\'{a}ndor, J.: On some new Wilker and Huygens type trigonometric-hyperbolic inequalities.  Proc. Jangjeon Math. Soc. \textbf{15}, 145--153 (2012)

\bibitem {s4}S\'{a}ndor, J.: Trigonometric and hyperbolic inequalities. Available via arXiv:1105.0859v1 [math.CA]

\bibitem {sh}Shafer,  R.E.: Problem E 1867.   Amer. Math. Monthly \textbf{73}, 309--310 (1966)

\bibitem {sh2}Shafer,  R.E.: On quadratic approximation. SIAM J. Numer. Anal. \textbf{11}, 447--460 (1974)

\bibitem {sh3}Shafer,  R.E.: Analytic inequalities obtained by quadratic approximation.  Univ. Beograd. Publ. Elektrotehn. Fak. Ser. Mat. Fiz. No. 577--598, 96--97 (1977)

\bibitem {sgmk} Shafer, R.E., Grinstein,  L.S., Marsh, D.C.B., Konhauser, J.D.E.: Problems and solutions: an inequality for the inverse tangent: E 1867. Amer. Math. Monthly \textbf{74}, 726--727 (1967)

\bibitem {sv2}Simi\'{c}, S., Vuorinen, M.: On quotients and differences of hypergeometric functions. J. Inequal. Appl. \textbf{2011}, Article 141, 10 pp.  (2011)

\bibitem {sv1}Simi\'{c}, S., Vuorinen, M.: Landen inequalities for zero-balanced hypergeometric functions.  Abstr. Appl. Anal. \textbf{2012}, Article ID 932061, 11 pp.  (2012)

\bibitem {st} Stolarsky,  K.B.: Generalizations of the logarithmic mean.  Math. Mag. \textbf{48},  87--92 (1975)

\bibitem{ta} Takeuchi,  S.: Generalized Jacobian elliptic functions and their application to bifurcation problems associated with $p$-Laplacian.  J. Math. Anal. Appl. \textbf{385}, 24--35 (2012)

\bibitem {to}Toader, Gh.: Some mean values related to the arithmetic-geometric mean.  J. Math. Anal. Appl. \textbf{218}, 358-368 (1998)

\bibitem {t} Trif, T.: Note on certain inequalities for means in two  variables.  JIPAM. J. Ineq. Pure Appl. Math. \textbf{6}, Article 43 (2005)

\bibitem {v} Vuorinen, M.: Hypergeometric functions in geometric function theory. In:  K. Srinivasa Rao, R. Jagannathan, G. Vanden Berghe, J. Van der Jeugt (eds.), Special functions and Differential equations, Proceedings of a workshop held at The Institute of Mathematical Sciences, Madras, India, January 13-24 (1997), Allied Publishers, pp. 119-126 (1998)

\bibitem {vz}Vuorinen, M., Zhang,  X.-H.: On exterior moduli of quadrilaterals and special functions. J. Fixed Point Theory Appl. (to appear).  Available via arXiv:1111.3812v2 [math.CA]

\bibitem {wc} Wang, M.-K., Chu, Y.-M.: Asymptotical bounds for complete elliptic integrals of the second kind. Available via arXiv:1209.0066v1 [math.CA]

\bibitem {wcj} Wang, M.-K., Chu, Y.-M., Jiang, Y.-P.: Ramanujan's cubic transformation inequalities for zero-balanced hypergeometric functions. Available via arXiv:1210.6126v1 [math.CA]

\bibitem {wcqj} Wang, M.-K., Chu, Y.-M., Qiu,  S.-L., Jiang, Y.-P.: Bounds for the perimeter of an ellipse.  J. Approx. Theory \textbf{164}, 928--937 (2012)

\bibitem {wqzc}Wang, G.-D.,  Qiu, S.-L.,  Zhang, X.-H., Chu, Y.-M.: Approximate convexity and concavity of generalized Gr\"otzsch ring function.  Appl. Math. J. Chinese Univ. Ser. B \textbf{21}, 203--206 (2006)

\bibitem {wzc1}Wang, G.-D.,  Zhang, X.-H., Chu, Y.-M.: Inequalities for the generalized elliptic integrals and modular functions.  J. Math. Anal. Appl. \textbf{331}, 1275--1283 (2007)

\bibitem {wzc3}Wang, G.-D.,  Zhang, X.-H., Chu, Y.-M.: A H\"older mean inequality for the Hersch-Pfluger distortion function. Sci. Sin. Math.  \textbf{40}, 783--786 (2010)

\bibitem {wzc2}Wang, G.-D.,  Zhang, X.-H., Chu, Y.-M.:  A power mean inequality for the Gr\"otzsch ring function. Math. Inequal. Appl. \textbf{14}, 833--837 (2011)

\bibitem {wzc4}Wang, G.-D.,  Zhang, X.-H., Chu, Y.-M.: A power mean inequality involving complete elliptic integrals.  Rocky Mountain J. Math. (to appear)

\bibitem  {wzj1}Wang, G.-D.,  Zhang, X.-H., Jiang, Y.-P.: Concavity with respect to H\"older means involving the generalized Gr\"otzsch function.  J. Math. Anal. Appl. \textbf{379}, 200--204 (2011)

\bibitem {wzj2}Wang, G.-D.,  Zhang, X.-H., Jiang, Y.-P.: H\"older concavity and inequalities for Jacobian elliptic functions.  Integral Transforms Spec. Funct. \textbf{23}, 337--345 (2012)

\bibitem  {wzqc}Wang, G.-D.,  Zhang, X.-H., Qiu, S.-L., Chu, Y.-M.: The bounds of the solutions to generalized modular equations.  J. Math. Anal. Appl. \textbf{321}, 589--594 (2006)

\bibitem {wi} Widder, D.V.: The Laplace Transform.  Princeton University Press (1941)

\bibitem {wr}  Wilker, J.B.: Problem E 3306.  Amer. Math. Monthly \textbf{96}, 55 (1989)

\bibitem {sjva} Wilker, J.B., Sumner, J.S., Jagers, A.A., Vowe, M., Anglesio, J.: Problems and solutions: solutions of elementary problems: E 3306.  Amer. Math. Monthly \textbf{98}, 264--267 (1991)

\bibitem  {ws} Williams,  J.P.: Solutions of advanced problems: a delightful inequality 5642.   Amer. Math. Monthly \textbf{76}, 1153--1154 (1969)

\bibitem {wu} Wu, S.-H.: Generalization and sharpness of the power means inequality and their applications.  J. Math. Anal. Appl. \textbf{312}, 637--652 (2005)

\bibitem {wb} Wu, S.-H., Baricz, \'A.: Generalizations of Mitrinovi\'c, Adamovi\'c and Lazarevi\'c's inequalities and their applications.   Publ. Math. Debrecen \textbf{75}, 447--458 (2009)

\bibitem {wd} Wu, S.-H., Debnath, L.: Jordan-type inequalities for differentiable functions and their applications.  Appl. Math. Lett. \textbf{21}, 803--809 (2008)

\bibitem {wd3} Wu, S.-H., Debnath, L.: Inequalities for differences of power means in two variables.  Anal. Math. \textbf{37}, 151--159 (2011)

\bibitem {ws1} Wu, S.-H., Srivastava, H.M.: A weighted and exponential generalization of Wilker's inequality and its applications.  Integral Transforms Spec. Funct. \textbf{18}, 529--535 (2007)

\bibitem {ws2} Wu, S.-H., Srivastava, H.M.: A further refinement of a Jordan type inequality and its application.  Appl. Math. Comput. \textbf{197}, 914--923 (2008)

\bibitem {YangSh}  Yang, S.-J.: Absolutely (completely) monotonic functions and Jordan-type inequalities.  Appl. Math. Lett. \textbf{25}, 571--574 (2012)

\bibitem {y} Yang, Z.-H.:  New sharp bounds for identric mean in terms of logarithmic mean and arithmetic mean.  J. Math. Inequal. \textbf{6}, 533--543 (2012)

\bibitem {ye}Yee, A. J.: Large computations (2010). Available at $http://www.\linebreak numberworld.org/nagisa\_runs/computations.html.$

\bibitem {yin}Yin, L.: Several inequalities for the volume of the unit ball in $\mathbb{R}^n$.  Bull. Malays. Math. Sci. Soc. (2) (to appear)

\bibitem{zwc3}Zhang, X.-H., Wang, G.-D., Chu,  Y.-M.: Some inequalities for the generalized Gr\"otzsch function.  Proc. Edinb. Math. Soc. (2) \textbf{51}, 265--272 (2008)

\bibitem {zwc}Zhang, X.-H., Wang, G.-D., Chu,  Y.-M.: Convexity with respect to H\"older mean involving zero-balanced hypergeometric functions.  J. Math. Anal. Appl. \textbf{353}, 256--259 (2009)

\bibitem {zwc2} Zhang, X.-H., Wang, G.-D., Chu,  Y.-M.: Remarks on generalized elliptic integrals.  Proc. Roy. Soc. Edinburgh Sect. A \textbf{139}, 417--426 (2009)

\bibitem {zwcq} Zhang, X.-H., Wang, G.-D., Chu,  Y.-M., Qiu, S.-L.: Monotonicity and inequalities for the generalized $\eta$-distortion function.  (Chinese) Chinese Ann. Math. Ser. A \textbf{28}, 183--190  (2007); translation in Chinese J. Contemp. Math. \textbf{28}, 141--148 (2007)

\bibitem {zcl}Zhao, T.-H., Chu,  Y.-M.,  Liu, B.-Y.: Some best possible inequalities concerning certain bivariate means. Available via arXiv:1210.4219v1 [math.CA]

\bibitem {zgq} Zhao, J.-L., Guo,  B.-N.,  Qi, F.: A refinement of a double inequality for the gamma function.  Publ. Math. Debrecen \textbf{50}, 1--10 (2011)

\bibitem {zqw} Zhou, L.-M., Qiu, S.-L., Wang, F.: Inequalities for the generalized elliptic integrals with respect to H\"{o}lder means.  J. Math. Anal. Appl.  \textbf{386}, 641--646 (2012)

\bibitem {z3} Zhu, L.: A general refinement of Jordan-type inequality.  Comput. Math. Appl. \textbf{55}, 2498--2505 (2008)

\bibitem {z7} Zhu, L.: On a quadratic estimate of Shafer.  J. Math. Inequal. \textbf{2}, 571--574 (2008)

\bibitem {z2} Zhu, L.: A general form of Jordan's inequalities and its applications.  Math. Inequal. Appl. \textbf{11}, 655--665 (2008)

\bibitem {z13} Zhu, L.: New inequalities of Shafer-Fink type for arc hyperbolic sine.  J. Inequal. Appl. \textbf{2008}, Article ID 368275, 5 pp. (2008)

\bibitem  {z5} Zhu, L.: Some new inequalities of the Huygens type. Comput. Math. Appl. \textbf{58}, 1180--1182 (2009)

\bibitem {z12} Zhu, L.: Some new Wilker-type inequalities for circular and hyperbolic functions. Abstr. Appl. Anal. \textbf{2009}, Article ID 485842, 9 pp. (2009)

\bibitem {z6} Zhu, L.: A source of inequalities for circular functions.  Comput. Math. Appl. \textbf{58}, 1998--2004 (2009)

\bibitem {z14} Zhu, L.: Generalized Lazarevi\'c's inequality and its applications -- Part II.  J. Inequal. Appl. \textbf{2009}, Article ID 379142, 4 pp.  (2009)

\bibitem {z15} Zhu, L.: Sharpening Redheffer-type inequalities for circular functions. Appl. Math. Lett. \textbf{22}, 743--748 (2009)

\bibitem {z4} Zhu, L.: Jordan type inequalities involving the Bessel and modified Bessel functions. Comput. Math. Appl. \textbf{59}, 724--736 (2010)

\bibitem {z10} Zhu, L.: A general form of Jordan-type double inequality for the generalized and normalized Bessel functions.  Appl. Math. Comput. \textbf{215},  3802--3810 (2010)

\bibitem {z1} Zhu, L.: An extended Jordan's inequality in exponential type.  Appl. Math. Lett. \textbf{24}, 1870--1873 (2011)

\bibitem {z16} Zhu, L.: Extension of Redheffer type inequalities to modified Bessel functions.  Appl. Math. Comput.  \textbf{ 217}, 8504--8506 (2011)

\bibitem  {zh} Zhu, L., Hua, J.-K.: Sharpening the Becker-Stark inequalities.   J. Inequal. Appl. \textbf{2010}, Article ID 931275, 4 pp. (2010)

\bibitem {zs} Zhu, L., Sun, J.-J.: Six new Redheffer-type inequalities for circular and hyperbolic functions.  Comput. Math. Appl. \textbf{56}, 522--529 (2008)


\end{thebibliography}
%

\end{document}